\documentclass[reqno, 10pt]{amsart}
\usepackage{a4wide}
\usepackage[greek, russian, german, french, italian, english]{babel}

\usepackage[T3, T1]{fontenc}
\usepackage[utf8]{inputenc}
\usepackage[italian]{varioref}

\usepackage{CJKutf8} 
\newcommand{\ZhSimplified}[1]{\begin{CJK}{UTF8}{gbsn}#1\end{CJK}}
\newcommand{\ZhTraditional}[1]{\protect\begin{CJK*}{UTF8}{bsmi}#1\end{CJK*}}

\usepackage{todonotes}

\usepackage{titlesec}
\usepackage{titletoc}

\renewcommand{\footnotesize}{\scriptsize}
\renewcommand{\thefootnote}{\alph{footnote}}

\usepackage[perpage]{footmisc}

\newcommand{\symbolfootnote}[1]{%
\let\oldthefootnote=\thefootnote%
\stepcounter{mpfootnote}%
\addtocounter{footnote}{-1}%
\renewcommand{\thefootnote}{\fnsymbol{mpfootnote}}%
\footnote{#1}%
\let\thefootnote=\oldthefootnote%
}

\usepackage{emptypage}
\usepackage{microtype} 
\usepackage{fancyhdr}
\usepackage{indentfirst}
\usepackage{changepage}

\usepackage{empheq} 

\usepackage{tikz-cd}
\usetikzlibrary{decorations.fractals}
\usetikzlibrary{lindenmayersystems} 

\usepackage{booktabs} 
\usepackage{tabularx} 
\usepackage{longtable}
\usepackage{graphicx} 
\usepackage{pict2e, curve2e}
\usepackage{listings}
\usepackage{enumitem}

\usepackage[all]{xy}

\usepackage{caption}
\usepackage{xcolor}
\definecolor{eggplant}{HTML}{800080}
\definecolor{mallard}{HTML}{008080}
\definecolor{duskycerulean}{HTML}{004080} 
\definecolor{dandelion}{HTML}{FDBC42}
\definecolor{pumpkin}{HTML}{FF7518}
\definecolor{vernalgreen}{HTML}{03E364}
\definecolor{magenta}{HTML}{E30382}
\definecolor{bluepoint}{HTML}{0382E3}

\usepackage{amsfonts}
\usepackage{amsmath}
\usepackage{amssymb}
\usepackage{amsaddr}
\usepackage{amsthm}
\usepackage{mathrsfs}
\usepackage[scr]{rsfso} 
\usepackage{esstixfrak}
\usepackage{eufrak}
\usepackage{fixmath}

\usepackage{scalerel}
\newlength\bshft
\def\pseudobold#1{\ThisStyle{\ooalign{$\SavedStyle#1$\cr%
  \kern-\bshft$\SavedStyle#1$\cr%
  \kern\bshft$\SavedStyle#1$}}}
  
\usepackage[version=4]{mhchem}

\usepackage[bbgreekl]{mathbbol}
\DeclareSymbolFontAlphabet{\mathbb}{AMSb}
\DeclareSymbolFontAlphabet{\mathbbl}{bbold}
                        
\usepackage{mathtools}
\usepackage{latexsym}
\usepackage{siunitx}
\usepackage{icomma}
\usepackage{slashed}
\usepackage{braket}
\usepackage{xfrac} 

\usepackage{bbding}

\usepackage{textcomp}

\titleformat{\section}[block]
{\large
\bfseries}
{\thesection.}{0.5em}{}

\titleformat{\subsection}[block]
{\normalsize
\bfseries}
{\thesubsection.}{0.5em}{}

\titleformat{\subsubsection}[block]
{\normalsize
\bfseries}
{\thesubsubsection.}{0.5em}{}

\contentsmargin{2.65em}

\titlecontents{chapter}
				[0pc]
				{\addvspace{.8pc}}
				{\bfseries\thecontentslabel{. }}
				{\hspace*{0em}}
				{\hfill\contentspage}
				[\addvspace{.3pc}]
				
\titlecontents{section}
				[0pc]
				{}
				{\thecontentslabel{. }}
				{\hspace*{0em}}
				{\titlerule*[0.4pc]{\thinspace.}\contentspage}

\titlecontents{subsection}
				[1.65pc]
				{}
				{\thecontentslabel{. }}
				{\hspace*{0em}}
				{\titlerule*[0.4pc]{\thinspace.}\contentspage}

\titlecontents{subsubsection}
				[3.95pc]
				{}
				{\thecontentslabel{. }}
				{\hspace*{0em}}
				{\titlerule*[0.4pc]{\thinspace.}\contentspage}

\newtheoremstyle{personal1}
  {\topsep}   
  {\topsep}   
  {\upshape}  
  {0pt}       
  {\itshape}  
  {.}         
  {5pt plus 1pt minus 1pt} 
  {\bfseries\thmname{#1}\thmnumber{ #2}\normalfont\thmnote{ \itshape(#3)}}
\theoremstyle{personal1}

\newtheoremstyle{personal2}
  {\topsep}   
  {\topsep}   
  {\itshape} 
  {0pt}      
  {\itshape} 
  {.}       
  {5pt plus 1pt minus 1pt} 
  {\bfseries\thmname{#1}\thmnumber{ #2}\normalfont\thmnote{ \itshape(#3)}} 
\theoremstyle{personal2}

\newtheorem{propositio}{Proposition}[section]

\newtheoremstyle{personal3}
  {\topsep}   
  {\topsep}   
  {\upshape}  
  {0pt}       
  {\itshape} 
  {.}         
  {5pt plus 1pt minus 1pt} 
  {\thmname{#1}\thmnumber{ \itshape#2}\thmnote{ (#3)}} 
\theoremstyle{personal3}

\newtheorem{exemplum}{Example}[section]
\newtheorem{marginalia}{Marginalia}[section]

\usepackage{float}

\graphicspath{figures}

\usepackage{url}

\usepackage[pdfauthor = {Edoardo{ }Niccolai},
	pdftitle = {scd},
		pdfsubject = {Mathematics{ }and{ }mathematical{ }physics}]{hyperref}
\hypersetup{colorlinks = true,
			citecolor = mallard,
			linkcolor = eggplant,
			urlcolor = duskycerulean}
			
\newcommand{\enumerationisinitium}{\begin{enumerate}[nolistsep, wide, label = \textnormal{($\mathnormal{\arabic{*}}$)}, ref = \textnormal{($\mathnormal{\arabic{*}}$)}]}
\newcommand{\enumerationisfinis}{\end{enumerate}}
\newcommand{\subenumerationisinitium}{\begin{enumerate}[nolistsep, wide, label = \textnormal{(\roman{*})}, ref = \textnormal{(\roman{*})}]}
\newcommand{\subenumerationisfinis}{\end{enumerate}}

\newenvironment{indent paragraph: 15pt}{%
  \par%
  \leftskip=15pt%
  \noindent\ignorespaces}{%
  \par}

\newenvironment{indent paragraph: 30pt}{%
  \par%
  \leftskip=30pt%
  \noindent\ignorespaces}{%
  \par}

\makeatletter
\newcommand{\xMapsto}[2][]{\ext@arrow 0599{\Mapstofill@}{#1}{#2}}
\def\Mapstofill@{\arrowfill@{\Mapstochar\Relbar}\Relbar\Rightarrow}
\makeatother

\DeclareSymbolFont{tipa}{T3}{cmr}{m}{n}
\DeclareMathAccent{\invertedbreve}{\mathalpha}{tipa}{16}

\DeclareMathOperator{\dAlembertian}{\square}

\DeclareMathOperator{\gradient}{grad}
\DeclareMathOperator{\Laplacian}{\bigtriangleup}
\DeclareMathOperator{\LittlewoodPaley}{\bigtriangleup}

\newcommand{\Besov}{B}

\newcommand{\Cl}{C\ell}

\newcommand{\Endomorphism}{\mathrm{end}}

\newcommand{\equivalent}{\stackrel{\textnormal{\tiny{eqv}}}{=\joinrel=}}

\newcommand{\F}{F}
\newcommand{\filtration}{\textcyrillic{\textit{Ф}}}

\newcommand{\idem}{\mathbb{I}}

\newcommand{\Lagrangian}{\mathscr{L}}

\newcommand{\Lebesgue}{L}
\newcommand{\length}{\ell}

\newcommand{\Liederivative}{\pounds}
\newcommand{\Lorentz}{\textcyrillic{\textit{Л}}}

\newcommand{\OrnsteinUhlenbeck}{\textgreek{\textit{Ϙ}}}

\newcommand{\Ric}{\mathrm{Ric}}

\newcommand{\Riemann}{R}

\newcommand{\rotatedc}{{\mathpalette\rotc\relax}}\newcommand{\rotc}[2]{\rotatebox[origin=c]{180}{$#1c$}}

\newcommand{\rotatedeta}{{\mathpalette\roteta\relax}}\newcommand{\roteta}[2]{\rotatebox[origin=c]{180}{$#1\eta$}}

\newcommand{\rotatedmcy}{{\mathpalette\rotmcy\relax}}\newcommand{\rotmcy}[2]{\rotatebox[origin=c]{180}{$#1\textcyrillic{\textit{м}}$}}

\newcommand{\rotatedPsi}{{\mathpalette\rotPsi\relax}}\newcommand{\rotPsi}[2]{\rotatebox[origin=c]{180}{$#1\Psi$}}

\newcommand{\rotatedupsilon}{{\mathpalette\rotupsilon\relax}}\newcommand{\rotupsilon}[2]{\rotatebox[origin=c]{180}{$#1\upsilon$}}
\newcommand{\rotatedvarphi}{{\mathpalette\rotvarphi\relax}}\newcommand{\rotvarphi}[2]{\rotatebox[origin=c]{180}{$#1\varphi$}}

\newcommand{\boundedvariation}{{\mathpalette\rotb\relax}}\newcommand{\rotb}[2]{\rotatebox[origin=c]{180}{$#1\textcyrillic{\textit{в}}$}}
\newcommand{\sadj}{\mathrm{s\text{-}adj}} 
\newcommand{\Schrodinger}{\ddot{S}}

\newcommand{\sezione}{\textgreek{\textit{ς}}}

\newcommand{\sigmaPauli}{\textgreek{\text{σ}}}

\newcommand{\Sobolev}{W}

\newcommand{\Tau}{\textgreek{Τ}}

\newcommand{\viz}{\stackrel{\textnormal{\tiny{viz}}}{=}}

\makeatletter

\newcommand*\sNeg[2][0mu]{\Neginternal{#1}{\snegslash}{#2}}

\newcommand*\Neginternal[3]{\mathpalette\Neg@{{#1}{#2}{#3}}}
\newcommand*\Neg@[2]{\Neg@@{#1}#2}
\newcommand*\Neg@@[4]{%
  \mathrel{\ooalign{%
    $\m@th#1#4$\cr
    \hidewidth$\m@th#3{#1}\mkern\muexpr#2*2$\hidewidth\cr
  }}%
}

\newcommand*\negslash[1]{\m@th#1\not\mathrel{\phantom{=}}}
\newcommand*\snegslash[1]{\rotatebox[origin=c]{60}{$\m@th#1-$}}
\newcommand*\ssnegslash[1]{\rotatebox[origin=c]{60}{$\m@th#1{\dabar@}\mkern-7mu{\dabar@}$}}
\newcommand*\sssnegslash[1]{\rotatebox[origin=c]{60}{$\m@th#1\dabar@$}}
\makeatother 

\makeatletter
\newcommand\footnoteref[1]{\protected@xdef\@thefnmark{\ref{#1}}\@footnotemark}
\makeatother

\usepackage{accents}
\newcommand*\underdot[1]{%
\underaccent{\dot}{#1}}

\addto{\captionsenglish}{\renewcommand{\refname}{Bibliography · List of Works Mentioned in the Article}
	}

\let\OLDthebibliography\thebibliography
\renewcommand\thebibliography[1]
{\OLDthebibliography{#1}
\setlength{\parskip}{0pt}
\setlength{\itemsep}{0pt plus 0.3ex}}

\title[]{Stochastic Covariant Derivatives \\ in a (Curved) Space-Time: \\ a Glimpse into the Fractoid Spaces}

\author{\small\href{https://edoardoniccolai.com}{Edoardo Niccolai}
	}

\begin{document}

\begin{abstract}
A study on the notion of covariant derivatives in flat and curved space-time via Itô–Wiener processes, when subjected to stochastic processes, is presented. Going into details, there is an analysis of the following topics:
(i) Besov space,
(ii) Schrödinger operators,
(iii) Klein–Gordon and Dirac equations,
(iv) Dirac operator via Clifford connection,
(v) semi-martingale and Stratonovich integral,
(vi) stochastic geodesics,
(vii) white noise on a (4+)D space-time $\mathfrak{H}$-geometry (with the Paley–Wiener integral), and
(viii) torsion of the covariant derivative. 
In the background stands the scale relativity theory, together with a sketch of the concept of fractoid spaces.

\vspace{2mm}

\noindent \textsc{Keywords}: $\sigma$-(sub)algebra(s), Besov space, $\mathscr{C}^2$ smooth vector fields, Clifford connection, covariant derivative of stochastic type in $\mathbb{R}^4_{1, 3}$ (Minkowski space-time) and in curved space-time (\textsc{gr}), Dirac equation(s), Dirac operator, forward and backward mean derivatives, fractal space-time, fractoid spaces, Itô–Wiener processes, Klein–Gordon equation(s), Littlewood–Paley operator, martingale, Newton–Nelson equation(s), Paley–Wiener integral, quantum-like fluctuations, random fluctuations, relativistic generalization of Nelson's stochastic mechanics, scale relativity, Schrödinger operators, Schwartz space, torsion.
\end{abstract}

\maketitle

\begingroup
\hypersetup{linktocpage}
\tableofcontents
\thispagestyle{empty}
\endgroup

\setcounter{secnumdepth}{0}  
\section{Thierry's Question}
\setcounter{secnumdepth}{3}

\begingroup
\footnotesize
Que devient la dérivée  covariante, que je note ici $D^*$, de Laurent  Nottale, obtenue en supposant un espace temps fractal, qui revient à un opérateur de type Schrödinger en $i\frac{d}{dt}$ + Laplacien dans un espace temps plat dans le cas non relativiste, dans le cas de la relativité général avec une métrique $g_{ij}$ qui fluctuerait de manière stochastique? \\
\indent — \textsc{T. Lehner}, via email to me

\endgroup

\vspace{2mm}

Patently, I do not offer a \textit{réponse univoque} to my sodalis Thierry, great and \textgreek{συμπαθητικός} friend, also because our starting perspectives are different (me with a mathematical priority, he with a physical priority), but a series of \textit{pièces de puzzle}, nay, a \ZhSimplified{七巧板} (qīqiǎobǎn), or \emph{tangram}, with a set of pieces, which, assembled together, can help and contribute to offering a réponse to his question.

\vspace{2mm}

\addcontentsline{toc}{section}{Ouverture}
\begin{center}
\textbf{\textsc{ouverture}}
\end{center}

\vspace{-2mm}

\section{Scale Relativity and Non-differentiability of a Continuous Fractal Space-Time}

\subsection{A Quick Summary}

To begin with, a few words about the scale relativity by L. Nottale \& collaborators \cite{Nottale "Scale Relativity Fractal Space-Time and Quantum Mechanics"} \cite{Nottale "Scale Relativity and Fractal Space-Time: Applications to Quantum Physics Cosmology and Chaotic Systems"} \cite{Nottale Celerier and T. Lehner Non-Abelian gauge field theory in scale relativity"} \cite{Nottale "Scale Relativity and Fractal Space-Time: Theory and Applications"} \cite{Nottale "Scale Relativity and Fractal Space-Time. A New Approach to Unifying Relativity and Quantum Mechanics"} \cite{Nottale and Lehner "Turbulence and Scale Relativity"}. This theory is made up of the following postulates.
\enumerationisinitium
\item The concept of space-time is \emph{continuous} but \emph{non-differentiable}, that is to say, \emph{fractal}. The space contemplated in the theory of scale relativity is thereby originated in the Mandelbrot set, 
\begin{equation}
	M_\mathbb{C} = \bigl\{c \in \mathbb{C} \mid \varphi_c^n(0) \sNeg[1mu]{\rightarrow} \infty \text{ as } n \to \infty\bigr\}.
\end{equation}
This set is related to a complex quadratic polynomial $\varphi_c$,
\begin{subequations}
	\begin{align}
	& 
	\label{align "Complex quadratic polynomial: map"}
	\mathbb{C} \xrightarrow{\varphi_c} \mathbb{C}, \\
	&
	\label{align "Complex quadratic polynomial: equation"}
	\varphi_c(z) = z^2 + c,
	\end{align}
\end{subequations}
under the map \eqref{align "Complex quadratic polynomial: map"} plus the equation \eqref{align "Complex quadratic polynomial: equation"}, hence a quadratic family $\varphi_c \colon z \mapsto z^2 + c$.

Consequently, 
\subenumerationisinitium
\item any covariant derivative is a math-construct of the non-differentiable and fractal geometry;
\item the baggage of coordinate transformations are continuous but can be non-differentiable.
\subenumerationisfinis
\item Let us try to find out what that means. 
\subenumerationisinitium
\item
\label{item "Curve and fractality"} 
A continuous curve line $\gamma_\mathrm{c}$, or rather, a geodesic, is non-differentiable when its length $\length(\gamma_\mathrm{c})$ is dependent on a scale variable $\varepsilon_\textsc{s}$, and it diverges, $\length\bigl({\varepsilon_\textsc{s}}_{(\gamma_\mathrm{c})}\bigr) \to \infty$, if $\varepsilon_\textsc{s} \to 0$, which constitutes its fractality; the scale divergence of continuous and almost nowhere-differentiable curves turns up as a counterpart to the extension of the fundamental theorem of calculus, of Lebesgueian memory, on the basis of which a curve, of finite length, is almost everywhere differentiable. For a continuous fractal $\alpha$-function, namely a scale-dependent function, $\varphi^\alpha(x) \text{ viz. } x(t, \delta t)$, all this is exemplified by two equality procedures—via nonstandard analysis, in  which $\delta t$ is replaced with $dt$—enunciating the variation of the position vector $x$ of a particle, along a geodesic, between $t - dt$ and $t$, and between $t$ and $t + dt$,
\begin{align}
	&
	\label{align "Variation of the position vector... 1"} 
	x(t + dt)dt - x(t, dt) = v_+ (x, t)dt + \sigma_\mathrm{s}^+ (t, dt)\left(\frac{dt}{\tau_\textsc{eb}}\right)^\frac{1}{\mathrm{D_f}}, \\
	&
	\label{align "Variation of the position vector... 2"}
	x(t, dt) - x(t - dt)dt = v_- (x, t) dt + \sigma_\mathrm{s}^- (t, dt)\left(\frac{dt}{\tau_\textsc{eb}}\right)^\frac{1}{\mathrm{D_f}},	
\end{align}
where $\sigma_\mathrm{s}^+$ and $\sigma_\mathrm{s}^-$ are \emph{stochastic variables} exhibiting finite \emph{velocity fluctuations}, $\tau_\textsc{eb} = \frac{\hbar}{mc^2}$ is the Einstein–de Broglie time, and $\mathrm{D_f}$ is the fractal dimension (or Hausdorff or Hausdorff–Besicovitch dimension) of the path. 

If one wants, Eqq. \eqref{align "Variation of the position vector... 1"} \eqref{align "Variation of the position vector... 2"} stand in a statistical explication, under diffusion coefficients.
\item The same rule of the previous point \ref{item "Curve and fractality"} applies to the notion of manifold, thus to a topological space. 
\item Regarding the fractal space-time, it entails the scale dependence of the reference frames, with internal scale variables: obviously, it will consist of non-differentiable but also differentiable parts (in any case, we are dealing with a \emph{continuum} which, historiquement et par définition, is differentiable); the non-differentiable part, the fractal one, will be characterized by \emph{stochastic fluctuations} (the two parts are later combined together). It is conjectured that fluctuations in the fractal space-time lead to $\delta{D_m}(x, t)$-fluctuations of the \emph{coefficient of diffusion} 
\begin{equation}
	D_m = \left(\frac{1}{2t}\right)\int^{+\infty}_{-\infty}\mathrm{\Delta}^2\mathscr{P}(\mathrm{\Delta}, t)d\mathrm{\Delta} = \left(\frac{1}{t}\right)\int^{+\infty}_{-\infty}\frac{\mathrm{\Delta}^2}{2}\mathscr{P}(\mathrm{\Delta}, t)d\mathrm{\Delta}, 
\end{equation}
where $\mathscr{P}(\mathrm{\Delta}, \tau)$ is the probability density function, and $d\mathrm{\Delta}$ is the displacement, or the length of path, from $\mathrm{\Delta}$ to $\mathrm{\Delta} + d\mathrm{\Delta}$, in an interval of time $t$. Cf. Eq. (14.20) in \cite{Niccolai "Notes in Pure Mathematics and Mathematical Structures in Physics"}; from such a context to the theory of Brownian motion, and interconnected Chapman–Kolmogorov equation, Fokker–Planck equation, etc., is a flash.
 
$\mathrm{N}\!\!\mathrm{B}$. There is also another fact that brings stochasticity in evidence: since the geodesics of fractal space-time are infinite, its interpretation must be stochastic.
\subenumerationisfinis
\item The scale relativity is conceived in such a way as to reject a \emph{discrete} space-time, as it is asserted in the \emph{poussière de Cantor} \cite[pp. 55-62]{Mandelbrot "Les objets fractals. Forme hasard et dimension"},\footnote{
	{} See, in this respect, the interesting notion of \emph{fractal lacunarity} \cite{Mandelbrot "Measures of fractal lacunarity: Minkowski content and alternatives"}.
	}
or in the fractal \textit{caillé}.
\enumerationisfinis

\subsection{The Head-scratcher of the Covariant Derivative}
\label{subsection "The Head-scratcher of the Covariant Derivative"}

The quandary of the covariant derivative in scale relativity is due to two conditions.

\subsubsection{Lack of Differentiability}

The first is that, where there is a lack of differentiability, owing to fractality, of space-time, one is witnessing a divergence that cannot be managed, except with stochastic tools—this is a clue of the \emph{mathematical skeleton} of the scale-dependent properties, which have a random bosom.\footnote{
	{} It is not a coincidence that the adjective \textgreek{στοχαστικός} means “able to hit”, “skilful in aiming”, “guessing”.
	}

\subsubsection{The Obstacle of the Two Derivatives, and the Arrangement with Sobolev and Besov Spaces}

The second condition is that \emph{two derivatives} arise, the classical one and the scale covariant derivative, the latter inserted in the fractal structure, once one accepts the hypothesis that the space-time is continuous, as the theory of scale relativity requires. 

That explains why one of the modes of enunciating \emph{fractional derivatives}, in favor of a good mathematical rigour, is to utilize  the standard Sobolev space \cite{Sobolev "Sur un theoreme d'analyse fonctionnelle"} \cite{Sobolev "Some Applications of Functional Analysis in Mathematical Physics"} 
\[
	\Sobolev^{k, 2}(Q) = \mathfrak{H}^k(Q), \Lebesgue^2(Q) = \mathfrak{H}^0(Q), \enspace k \in \mathbb{Z}_*, \text{ or } k \in \mathbb{N}, Q \subset \mathbb{R}^n, 
\]	
easily tractable; for $p = 2$, the Sobolev space $\Sobolev^{k, p}(Q)$ is a Hilbert space.

An analogous goal, toward the \emph{fractal Laplacian}, is achieved with the Besov space \cite{Besov Il'In Nikol'skii "Integral Representations of Functions and Imbedding Theorems I"} \cite{Besov Il'In Nikol'skii "Integral Representations of Functions and Imbedding Theorems II"}\footnote{
	{} For a detailed bibliography on Besov spaces, see \cite[pp. 898-901]{Sawano "Theory of Besov Spaces"}.
	}
\[
	\Besov^s_{p, \infty}(\mathbb{R}^2), \enspace s \subset \mathbb{R}^n, 0 < p \leqslant \infty, 
\]
which generalizes the Sobolev space, when the $\Besov$-space is a collection of every regular tempered (or Schwartz) distribution $\upsilon_\text{t} \in \mathcal{D}^*(\mathbb{R}^2)$, $\mathcal{D}^* \equivalent \mathcal{S}^*_\mathfrak{c}$ being the Schwartz space \cite{Grothendieck "Sur les espaces (F) et (DF)"} \cite{Dieudonne Schwartz "La dualite dans les espaces (F) et (LF)"} \cite{Schwartz "Produits tensoriels topologiques d'espaces vectoriels topologiques. Espaces vectoriels topologiques nucleaires"} of distributions, such that the norms
\begin{subequations}
\begin{align}
	& 
	\label{align "Besov space norm 1"}
	\left\|\,\rotatedc \mid \Lebesgue^p(\mathbb{R}^2)\right\| + \sup_{0 < |\upsilon_\text{t}| \leqslant 1} |\upsilon_\text{t}|^{-s} \left\|\LittlewoodPaley^k_{\upsilon_\text{t}}\rotatedc \mid \Lebesgue^p(\mathbb{R}^2)\right\|, \enspace \tfrac{2(1 - p)}{p} < s < k \in \mathbb{N}, \\
	&
	\label{align "Besov space norm 2"} 
	\left\|\,\rotatedc \mid \Lebesgue^p(\mathbb{R}^2)\right\| +  \left(\int^1_0|\upsilon_\text{t}|^{-sq}\left\|\LittlewoodPaley^k_{\upsilon_\text{t}}\rotatedc \mid \Lebesgue^p(\mathbb{R}^2)\right\|^q\frac{d\upsilon_\text{t}}{|\upsilon_\text{t}|^2}\right)^\frac{1}{q}, \enspace 0 < q < \infty,
\end{align}
\end{subequations}
are finite, id est Eqq. \eqref{align "Besov space norm 1"} \eqref{align "Besov space norm 2"} $< \infty$, where $\rotatedc$ is a continuous monotonically increasing function, and $\Lebesgue^p$ the Lebesgue space. The operator $\LittlewoodPaley$ can be recognized as the Littlewood–Paley operator \cite{Littlewood Paley "Theorems on Fourier Series and Power Series"} \cite{Littlewood Paley "Theorems on Fourier Series and Power Series II"} \cite{Littlewood Paley "Theorems on Fourier Series and Power Series III"}.

\subsubsection{Do Not Forget the Discretum}

I disagree with maintaining the continuity of space as the womb of space-time. Actually, the continuum seems to be an \emph{approximation of the discretum}, cf. \cite[Margo 9.2.1]{Niccolai "Notes in Pure Mathematics and Mathematical Structures in Physics"}. This is the \emph{original sin} inherent in the theory of scale relativity, as far as I am able to judge. 

\vspace{2mm}

\addcontentsline{toc}{section}{Parts of the Answer: a Fab-Tangram}
\begin{center}
\textbf{\textsc{parts of the answer: a fab-tangram}}
\end{center}

\section{Some Tips. Let us Do a Refresh: Schrödinger Operators}
\label{section "Some Tips. Let us Do a Refresh: Schrödinger Operators"}

Let us first give some definitions, just to cross some T's and dot some I's.

\subsection{Linear Schrödinger Operator for a 1-parameter Unitary Group}

The simple equality 
\begin{equation}
\label{equation "The Schrödinger operator for a 1-parameter unitary group"}
	\Schrodinger_t \viz U_t = e^{it\Laplacian}
\end{equation}
establishes the Schrödinger operator, whose notation here is a letter $S$ with diaeresis, for a 1-parameter unitary group. The symbol $\Laplacian$ designates the Laplacian.

For an equality 
\[
	\Schrodinger_t\varphi = \psi(\hat{\mathrm{x}}, t), 
\]	
the linear Schrödinger operator \eqref{equation "The Schrödinger operator for a 1-parameter unitary group"} is determined by
\begin{equation}
	\Schrodinger_t\varphi(\hat{\mathrm{x}}) = \left(\frac{1}{4\pi it}\right)^\frac{n}{2}\int e^{i\frac{|\hat{\mathrm{x}} - \dot{\hat{\mathrm{x}}}|^2}{4t}}\varphi \left(\dot{\hat{\mathrm{x}}}\right)d\dot{\hat{\mathrm{x}}}, \enspace \hat{\mathrm{x}} \in \mathbb{R}^n.
\end{equation}

\subsection[Schrödinger Operator on Continuous $\Lebesgue^2(\mathbb{R}^n)$-Space Structures]{Schrödinger Operator on Continuous $\protect\pseudobold{\Lebesgue^2(\mathbb{R}^n)}$-Space Structures}

The Schrödinger operator
\begin{equation}
\label{equation "Schrödinger operator on continuous..."}
	\Schrodinger = - \frac{\hbar^2}{2}\Laplacian + \rotatedupsilon[\mathbb{R}].
\end{equation}
is a linear partial differential operator in the Hilbert space—and in fact $\mathfrak{H} \viz \Lebesgue^2(\mathbb{R}^n)$; in Eq. \eqref{equation "Schrödinger operator on continuous..."} the Laplacian is $n$-dimensional, and $\rotatedupsilon$ is a smooth $\mathbb{R}$-potential.
 
\subsection[Covariant Schrödinger Operator on Riemannian $n$-Manifolds (Groupal Algebra)]{Covariant Schrödinger Operator on Riemannian $\protect\pseudobold{n}$-Manifolds (Groupal Algebra): $\protect\pseudobold{\Schrodinger^\nabla_\textcyrillic{\textit{э}}}$ in $\protect\pseudobold{\Gamma_{\sezione}^{(\Lebesgue^2)}}$}

We define a covariant Schrödinger operator $\Schrodinger^\nabla_\textcyrillic{\textit{э}}$, where $\nabla$ is a (metric) covariant derivative on a (metric) vector bundle $\mathring{\mathcal{E}}$ over $\mathcal{M}$, i.e. $\mathring{\mathcal{E}} \to \mathcal{M}$, and $_{\sadj}(\textcyrillic{\textit{э}}) \viz \textcyrillic{\textit{э}}$ is a self-adjoint endomorphism, on a smooth Riemannian pair $(\mathcal{M}, g)$ in the Hilbert space $\Lebesgue^2$ of square-integrable sections, denoted by 
\[
	\Gamma_{\sezione}^{(\Lebesgue^2)}\left(\mathcal{M}, \mathring{\mathcal{E}}\right), 
\]
so to get
\begin{equation}
\left(e^{-t\Schrodinger^\nabla_\textcyrillic{\textit{э}}}(x)\right)_{t \geqslant 0}	 \subset \Gamma_{\sezione}^{(\Lebesgue^2)}\left(\mathcal{M}, \mathring{\mathcal{E}}\right), \enspace t > 0, x \in \mathcal{M}.
\end{equation}
Let $d\bbmu$ be the Riemannian volume element ($\bbmu$ is for a Borel measure); for a function 
\[
	\varphi_\bbmu \in \Gamma_{\sezione}^{(\Lebesgue^2)}\left(\mathcal{M}, \mathring{\mathcal{E}}\right), 
\]
one has
\begin{equation}
	e^{-t\Schrodinger^\nabla_\textcyrillic{\textit{э}}}\varphi(x) = \int_\mathcal{M}e^{-t\Schrodinger^\nabla_\textcyrillic{\textit{э}}}(x, y)\varphi_\bbmu(y)d\bbmu(y).
\end{equation}
Here the Schrödinger operator $\Schrodinger^\nabla_\textcyrillic{\textit{э}}$ is but a \emph{covariant Schrödinger bundle}, that is, 
\[
	\left(\mathring{\mathcal{E}}, \nabla, \textcyrillic{\textit{э}}\right) \xrightarrow{(\cdot)} \mathcal{M},
\]
bearing in mind that the map
\begin{equation}
\label{equation "Endomorphism"}
	\textcyrillic{\textit{э}} \colon \mathcal{M} \to \text{end}\left(\mathring{\mathcal{E}}\right)
\end{equation}
is Borel $\bbmu$-measurable, having a linear self-adjoint map $\textcyrillic{\textit{э}}(x) \colon \mathring{\mathcal{E}}_x \to \mathring{\mathcal{E}}_x$, for each $x \in \mathcal{M}$. The self-adjoint endomorphism $_{\sadj}(\textcyrillic{\textit{э}})$ in \eqref{equation "Endomorphism"} counts as a potential on $\mathring{\mathcal{E}} \to \mathcal{M}$.

\subsection[Random Schrödinger Operator on Discrete $\Lebesgue^2(\mathbb{Z}^n)$-Space Structures]{Random Schrödinger Operator on Discrete $\protect\pseudobold{\Lebesgue^2(\mathbb{Z}^n)}$-Space Structures}

If one needs to use a random Schrödinger operator can quietly go from a $\Lebesgue^2(\mathbb{R}^n)$-space to $\Lebesgue^2(\mathbb{Z}^n)$-space, in such a way that the $\Laplacian$-operator from continuous on the $\mathbb{R}$-field becomes discrete on the $\mathbb{Z}$-field, so as to have this distinctness:
\begin{equation}
	-\Laplacian_\mathbb{Z} = \sum^n_{\lambda = 1}\Bigl({^2\varphi}(x)_\mathbbl{z} - \varphi(x - e_\lambda)_\mathbbl{z} - \varphi(x + e_\lambda)_\mathbbl{z}\Bigr).
\end{equation}

\section{Some Answers—Let Us Just Cut to the Chase: Covariant Derivatives}

From this Section, until the end of the articles, I will sketch some solutions to the opening question.

\subsection[Covariant Derivative of Stochastic Type in $\mathbb{R}^4_{1, 3}$ (c.-à-d. on a Flat Lorentz–Minkowski Space-Time) via Itô–Wiener Processes]{Covariant Derivative of Stochastic Type in $\protect\pseudobold{\mathbb{R}^4_{1, 3}}$ (c.-à-d. on a Flat Lorentz–Minkowski Space-Time) via Itô–Wiener Processes}
\label{subsection "Covariant Derivative of Stochastic Type (on a Flat Lorentz–Minkowski Space-Time) in Minkowski Space-Time via Itô–Wiener Processes"}
\markright{Stochastic Covariant Derivative in $\mathbb{R}^4_{1, 3}$ (on a Flat Lorentz–Minkowski Space-Time) via Itô–Wiener Processes}

Here we look for the covariant derivative, within the stochastic realm, in Minkowski/Lorentz–Minkowski space-time $\mathbb{M}^4 = \mathbb{R}^4_{1, 3}$. 

The first thing to do is recover Itô's formula, for $\mathbb{M}^4 \viz \mathcal{M}^4$, together with the Wiener processes, for $\mathbb{R}^4_{1, 3}$.

We choose
\begin{equation}
\label{equation "Itô diffusion-type process 1"}
	\rotatedvarphi_t = \rotatedvarphi_0 + \int^t_0 \beta_sds + \textcyrillic{\textit{в}}_\mathbbl{R}\mathsf{W}_t
\end{equation}
as an \emph{Itô diffusion-type process}, where 

$\rotatedvarphi_t$ is a stochastic process,
 
$\beta_s$ is a process almost surely with bounded variation of some path,  

$\textcyrillic{\textit{в}}_\mathbbl{R} > 0$ is a real constant, 

$\mathsf{W}_t$ is the Wiener stochastic process, which is almost surely continuous in $t$, and square-integrable martingale regarding a non-decreasing family $\mathscr{U}_t$, $t \in [0, \infty)$ of $\sigma$-subalgebras of the $\sigma$-algebra $\mathscr{B}$. 

\vspace{2mm}

\begin{marginalia}
A martingale is a stochastic process governed by a sequence of \emph{random fluctuations}.

Consider that in my formalism—cf. \cite[Sec. 12.4.3.2, and Definition 16.1.9]{Niccolai "Notes in Pure Mathematics and Mathematical Structures in Physics"}—the triple $(\invertedbreve{\Omega}, \mathscr{B}_\sigma, \bbmu)$ denotes a \emph{probability space} with a Borel $\sigma$-algebra on $\invertedbreve{\Omega}$.
\end{marginalia}	

\vspace{2mm}

Let 

$\tau$ be a proper time (an invariant parameter), 

$\textcyrillic{\textit{Э}} \text{ viz. } \textcyrillic{\textit{Э}}(\mathscr{U}^\rotatedvarphi_\tau)$ be the conditional expectation on $(\invertedbreve{\Omega}, \mathscr{B}_\sigma, \bbmu)$ of $\rotatedvarphi_{(t)}$ concerning the $\sigma$-algebra generated by some Borel sets in a $n$-dimensional $\mathbb{R}$-field, with the map $\upsilon \colon \invertedbreve{\Omega} \to \mathbb{R}^n$ (that is why the expectation covers the $\mathscr{U}_t$),

$\textit{\textgreek{Ν}}$ be the present state,\footnote{
	{} It is not a La. letter. It is a capital Gr. letter, from the adverb \textgreek{νῦν}, “now”.
	}
viz. the now, of $\rotatedvarphi_\tau$, scilicet the present $\sigma$-algebra for $\rotatedvarphi_\tau$.
 
Suppose that $\rotatedvarphi_\tau$ (the stochastic process) has values in a specific Riemannian manifold, in compliance with the map $\rotatedvarphi_\tau \colon \invertedbreve{\Omega} \to \mathcal{M}$. 

Then we can write, à la Dohrn–Guerra–Ruggiero \cite{Dohrn Guerra Ruggiero "Spinning Particles and Relativistic Particles in the Framework of Nelson's Stochastic Mechanics"} \cite{Guerra and Ruggiero "A Note on Relativistic Markov Processes"},\footnote{
	{} See also A.B. Cruzeiro and J.-C. Zambrini \cite[sec. 5]{Cruzeiro and Zambrini "Feynman's Functional Calculus and Stochastic Calculus of Variations"}.
	} 
the \emph{relativistic forward and backward mean derivative of stochastic type}, indicated with $D^+\rotatedvarphi_\tau$ and $D^-\rotatedvarphi_\tau$, respectively, for a flat pseudo-Euclidean (Minkowskian-like) Lorentzian space-time:
\begin{align}
\label{align "Forward mean derivative"}
	D^+\rotatedvarphi_\tau = & \lim_{\mathrm{\Delta}_\tau \downarrow 0}\textcyrillic{\textit{Э}}\left\{\frac{\rotatedvarphi(\tau + \mathrm{\Delta}\tau) - \rotatedvarphi_\tau}{\mathrm{\Delta}\tau} \mathrel{\bigg|} {\textit{\textgreek{Ν}}}^\rotatedvarphi_t\Bigl(\rotatedvarphi(\tau + \mathrm{\Delta}\tau) - \rotatedvarphi_\tau\Bigr)^2 \leqslant 0\right\} \notag \\
	& + \lim_{\mathrm{\Delta}_\tau\downarrow 0}\textcyrillic{\textit{Э}}\left\{\frac{\rotatedvarphi_\tau - \rotatedvarphi(\tau - \mathrm{\Delta}\tau)}{\mathrm{\Delta}\tau} \mathrel{\bigg|}{\textit{\textgreek{Ν}}}^{\rotatedvarphi}_t\Bigl(\rotatedvarphi_\tau - \rotatedvarphi(\tau - \mathrm{\Delta}\tau)\Bigr)^2 \geqslant 0\right\}
\end{align}
and
\begin{align}
\label{align "Backward mean derivative"}
	D^-\rotatedvarphi_\tau = & 
	\lim_{\mathrm{\Delta}_\tau \downarrow 0}\textcyrillic{\textit{Э}}\left\{\frac{\rotatedvarphi_\tau - \rotatedvarphi(\tau - \mathrm{\Delta}\tau)}{\mathrm{\Delta}\tau} \mathrel{\bigg|} {\textit{\textgreek{Ν}}}^{\rotatedvarphi}_t\Bigl(\rotatedvarphi_\tau - \rotatedvarphi(\tau - \mathrm{\Delta}\tau)\Bigr)^2 \leqslant 0\right\} \notag \\
	& + \lim_{\mathrm{\Delta}_\tau \downarrow 0}\textcyrillic{\textit{Э}}\left\{\frac{\rotatedvarphi(\tau + \mathrm{\Delta}\tau) - \rotatedvarphi_\tau}{\mathrm{\Delta}\tau} \mathrel{\bigg|} {\textit{\textgreek{Ν}}}^{\rotatedvarphi}_t\Bigl(\rotatedvarphi (\tau + \mathrm{\Delta}\tau) - \rotatedvarphi_\tau\Bigr)^2 \geqslant 0\right\},
\end{align}
marking with $\mathrm{\Delta}\tau$ the relativistic displacements a/o increments of $\rotatedvarphi_\tau$.

The relativistic forward mean derivative $D^+\rotatedvarphi_\tau$ \eqref{align "Forward mean derivative"} and the relativistic backward mean derivative $D^-\rotatedvarphi_\tau$ \eqref{align "Backward mean derivative"} are \emph{covariant} under the Lorentz transformations of the reference systems given by the tetrads. 
 
It will of course be useful to note that one has $D^+\rotatedvarphi_\tau = \rotatedPsi^+\Bigl(\tau, \rotatedvarphi_\tau\Bigr)$, for a  $\mathscr{C}^2$ smooth vector field
\begin{align}
	\rotatedPsi^+(\tau, x) = & \lim_{\mathrm{\Delta}_\tau \downarrow 0}\textcyrillic{\textit{Э}}\left\{\frac{\rotatedvarphi(\tau + \mathrm{\Delta}\tau) - \rotatedvarphi_\tau}{\mathrm{\Delta}\tau} \mathrel{\bigg|} \rotatedvarphi_\tau = x\Bigl(\rotatedvarphi(\tau + \mathrm{\Delta}\tau) - \rotatedvarphi_\tau\Bigr)^2 \leqslant 0\right\} \notag \\
	& + \lim_{\mathrm{\Delta}_\tau\downarrow 0}\textcyrillic{\textit{Э}}\left\{\frac{\rotatedvarphi_\tau - \rotatedvarphi(\tau - \mathrm{\Delta}\tau)}{\mathrm{\Delta}\tau} \mathrel{\bigg|}\rotatedvarphi_\tau = x\Bigl(\rotatedvarphi_\tau - \rotatedvarphi(\tau - \mathrm{\Delta}\tau)\Bigr)^2 \geqslant 0\right\},
\end{align}
and $D^-\rotatedvarphi_\tau = \rotatedPsi^-\Bigl(\tau, \rotatedvarphi_\tau\Bigr)$, for a $\mathscr{C}^2$ smooth vector field
\begin{align}
	\rotatedPsi^-(\tau, x) = & 
	\lim_{\mathrm{\Delta}_\tau \downarrow 0}\textcyrillic{\textit{Э}}\left\{\frac{\rotatedvarphi_\tau - \rotatedvarphi(\tau - \mathrm{\Delta}\tau)}{\mathrm{\Delta}\tau} \mathrel{\bigg|} \rotatedvarphi_\tau = x \Bigl(\rotatedvarphi_\tau - \rotatedvarphi(\tau - \mathrm{\Delta}\tau)\Bigr)^2 \leqslant 0\right\} \notag \\
	& + \lim_{\mathrm{\Delta}_\tau \downarrow 0}\textcyrillic{\textit{Э}}\left\{\frac{\rotatedvarphi(\tau + \mathrm{\Delta}\tau) - \rotatedvarphi_\tau}{\mathrm{\Delta}\tau} \mathrel{\bigg|} \rotatedvarphi_\tau = x\Bigl(\rotatedvarphi (\tau + \mathrm{\Delta}\tau) - \rotatedvarphi_\tau\Bigr)^2 \geqslant 0\right\}.
\end{align}

\subsection{Covariant Derivative of Stochastic Type in Curved Space-Time (c.-à-d. on a Lorentz Hyperbolic Manifold) via Itô–Wiener Processes}
\label{subsection "Covariant Derivative of Stochastic Type in Curved Space-Time (on a Lorentz Hyperbolic Manifold) via Itô–Wiener Processes"}
\markright{Stochastic Covariant Derivative in Curved Space-Time (on a Lorentz Hyperbolic Manifold) via Itô–Wiener Processes}

We analyze the context of general relativity (\textsc{gr}). Let $\rotatedPsi(\tau, m)$ be a vector field on a Lorentz 4-manifold, symbolized by $\mathbb{L}^4$, with a metric signatures $^{(1, 3)^-} \text{ viz. } (-, +, +, +)$, or, which is the same, to use the usual Riemannian notation, on a 4-manifold, symbolized by $\mathcal{M}^4$, of type $\mathscr{C}^2$ smooth, considering within this scheme a Lorentzian orthonormal frame in the tangent space $\mathcal{T}_x\mathcal{M}^4$, $x \in \mathcal{M}^4$. 

Let us say that $\Gamma_{\tau, s}$ is an operator of parallel translation—derived from the Levi-Civita-like connection—on the Lorentz bundle $\Lorentz(\mathcal{M}^4)$, along a stochastic Itô-process \cite{Ito "The Brownian motion and tensor fields on Riemannian manifold"} \cite{Ito "Stochastic parallel displacement"},\footnote{
	{} See the finishing touches in \cite{Dohrn and Guerra "Geodesic correction to stochastic parallel displacement of tensors"}.
	} 
from a \emph{random point} $\rotatedvarphi_s$ to another \emph{random point} $\rotatedvarphi_\tau$. Therefore the \emph{displacements a/o deviations of the geodesic} are taken into account. 
 
For the general relativity, we define the \emph{covariant relativistic mean derivatives of stochastic type}, 
\[
 	D^+\rotatedPsi(\tau, \rotatedvarphi_\tau) \text{ and } D^-\rotatedPsi(\tau, \rotatedvarphi_\tau),  
\]
over a $\mathbb{L}^4$-manifold, on the guideline of these equalities:
\begin{equation}
	D\rotatedPsi(\tau, \rotatedvarphi_\tau) = \lim_{\mathrm{\Delta}\tau \to + 0} \textcyrillic{\textit{Э}}^\rotatedvarphi_\tau\left\{\frac{\Gamma_{\tau, \tau + \mathrm{\Delta}\tau}\rotatedPsi\Bigl(\tau + \mathrm{\Delta}\tau, \rotatedvarphi(\tau + \mathrm{\Delta}\tau)\Bigr) - \rotatedPsi(\tau, \rotatedvarphi_\tau)\Bigr)}{\mathrm{\Delta}\tau}\right\},
\end{equation}
and
\begin{equation}
	D^*\rotatedPsi(\tau, \rotatedvarphi_\tau) = \lim_{\mathrm{\Delta}\tau \to + 0} \textcyrillic{\textit{Э}}^\rotatedvarphi_\tau\left\{\frac{\rotatedPsi(\tau, \rotatedvarphi_\tau) - \Gamma_{\tau, \tau - \mathrm{\Delta}\tau}\rotatedPsi\Bigl(\tau - \mathrm{\Delta}\tau, \rotatedvarphi(\tau - \mathrm{\Delta}\tau)\Bigr)}{\mathrm{\Delta}\tau}\right\},
\end{equation}
after specifying the expressions 
\begin{align}
	D\rotatedPsi = \frac{\partial\rotatedPsi}{\partial_\tau} + \nabla_\beta \rotatedPsi + \left(\frac{1}{2}\right)\nabla^2\rotatedPsi = \frac{\partial\rotatedPsi}{\partial_\tau} +\nabla_\beta \rotatedPsi + \left(\frac{\textcyrillic{\textit{в}}_\mathbbl{R}^2}{2}\right)\nabla^2\rotatedPsi, \\ 
	D^*\rotatedPsi = \frac{\partial\rotatedPsi}{\partial_\tau} + \nabla_{\beta *}\rotatedPsi - \left(\frac{1}{2}\right)\nabla^2\rotatedPsi = \frac{\partial\rotatedPsi}{\partial_\tau} + \nabla_{\beta *}\rotatedPsi - \left(\frac{\textcyrillic{\textit{в}}_\mathbbl{R}^2}{2}\right)\nabla^2\rotatedPsi,
\end{align}
where $\nabla$ is the covariant derivative of the Levi-Civita connection, $\nabla^2$ is the Laplace–Beltrami operator, and 
\begin{equation}
\label{equation "Real constant (squared) and the h-bar"}
	\frac{\textcyrillic{\textit{в}}_\mathbbl{R}^2}{2} \equivalent \frac{\hbar}{2m}. 
\end{equation}	
The reduced Planck constant, $\hbar = \frac{h}{2\pi}$, pops up on the assignment of the value 
\begin{equation}
\label{equation "Real constant and the h-bar"}
	\textcyrillic{\textit{в}}_\mathbbl{R} = \sqrt{\left(\frac{\hbar}{m}\right) \equivalent 1} 
\end{equation}
to the number $\textcyrillic{\textit{в}}_\mathbbl{R}$, embracing the Itô diffusion-type process, subordinated to the statement 
\begin{equation}
	\rotatedvarphi_{\left[\mathcal{M}^4\right]} = \int^\tau_0 \beta_sds + \textcyrillic{\textit{в}}_\mathbbl{R}\mathsf{W}_\tau,
\end{equation}
to be compared with Eq. \eqref{equation "Itô diffusion-type process 1"}.

The aforementioned decomposition is \emph{covariant} referring to Lorentz transformations in $\mathcal{T}_x\mathcal{M}^4$ for \eqref{align "Forward mean derivative"} and \eqref{align "Backward mean derivative"} bur not for $D^\rotatedvarphi_\tau$ and $D^*\rotatedvarphi_\tau$. There is however another covariance. Let $\mathcal{X}$ be a frame under which the time-like component and the space-like $3\mathrm{D}$ components hold, so 
\begin{align}
	& D^+\rotatedvarphi_\tau = \mathcal{X}^0_+(\tau, \rotatedvarphi_\tau), \\
	& D^-\rotatedvarphi_\tau = \mathcal{X}^0_-(\tau, \rotatedvarphi_\tau).
\end{align} 

Let $w^{\rotatedvarphi_\tau}$ be a vector field. Putting
\begin{align}
	& \tilde{w}^\rotatedvarphi_1(\tau, x) = \frac{1}{2}\Bigl(\mathcal{X}^0_+(\tau, x) + \mathcal{X}^0_-(\tau, x)\Bigr), \\
	& \tilde{w}^\rotatedvarphi_2(\tau, x) = \frac{1}{2}\Bigl(\mathcal{X}^0_+(\tau, x) - \mathcal{X}^0_-(\tau, x)\Bigr). 
\end{align}
We set 
\begin{equation}
	\tilde{w}^\rotatedvarphi_1(\tau, \rotatedvarphi_\tau) = \tilde{D}_\textsc{s}\rotatedvarphi_\tau 
\end{equation}	
as the \textsc{gr}-relativistic current 4-velocity of $\rotatedvarphi_\tau$,\footnote{
	\label{footnote "velocity and stochastic fluctuation"}
	When we talk about velocity, a \emph{stochastic fluctuation} takes place here.
	} 
where 
\begin{equation} 
	\tilde{D}_\textsc{s} = \frac{1}{2}(D^+ + D^-),
\end{equation}
is the \textsc{gr}-relativistic symmetric (\textsc{s}) mean derivative, and
\begin{equation}
	\tilde{w}^\rotatedvarphi_2(\tau, \rotatedvarphi_\tau) = \tilde{D}_\textsc{a}\rotatedvarphi_\tau 
\end{equation}
as the \textsc{gr}-relativistic velocity (cf. footnote \ref{footnote "velocity and stochastic fluctuation"} on p. \pageref{footnote "velocity and stochastic fluctuation"}) of osmotic determination of $\rotatedvarphi_\tau$, where
\begin{equation}
	\tilde{D}_\textsc{a} = \frac{1}{2}(D^+ - D^-)
\end{equation}
is the \textsc{gr}-relativistic antisymmetric (\textsc{a}) mean derivative.

The current vector $w^{\rotatedvarphi_\tau}_1$ is clearly covariant, since
\begin{equation}
	w^{\rotatedvarphi_\tau}_1 = \tilde{w}^{\rotatedvarphi_\tau}_1
\end{equation} 
(they have the the same decomposition-coordinate); whilst about the osmotic vector, it is decreed that
\begin{align}
	& 	w^{\rotatedvarphi_\tau}_2 = \left(D_\textsc{a}^\rotatedvarphi\rotatedvarphi^0_\tau, D_\textsc{a}^\rotatedvarphi\tilde{\rotatedvarphi}_\tau\right), \\
	& \tilde{w}^{\rotatedvarphi_\tau}_2 = \left(D_\textsc{a}^\rotatedvarphi\rotatedvarphi^0_\tau - D_\textsc{a}^\rotatedvarphi\tilde{\rotatedvarphi}_\tau\right).
\end{align}

The vector $\frac{1}{2}(D^+D^- + D^- D^+)\rotatedvarphi_\tau$ of (the stochastic process) $\rotatedvarphi_\tau$, is called the \emph{4-acceleration}. Under the Itô formula in $\mathcal{M}^4$ emerges that
\begin{equation}
	\frac{1}{2}(DD^* + D^*D)\rotatedvarphi_\tau = \frac{1}{2}(D^+D^- + D^- D^+)\rotatedvarphi_\tau,
\end{equation}
which is proven by $D^+D^-\rotatedvarphi_\tau = \bigl(D^*_\rotatedvarphi D_\rotatedvarphi\rotatedvarphi_\tau^0, D_\rotatedvarphi D^*_\rotatedvarphi\tilde{\rotatedvarphi}_\tau\bigr)$ and $D^-D^+\rotatedvarphi_\tau = \bigl(D_\rotatedvarphi D^*_\rotatedvarphi\rotatedvarphi_\tau^0, D^*_\rotatedvarphi D_\rotatedvarphi\tilde{\rotatedvarphi}_\tau\bigr)$.

\subsection{A Step Back}

Let us proceed with an additional conceptualization. The 4-vector acceleration $\frac{1}{2}(DD^* + D^*D)\rotatedvarphi_\tau$ has a precise stochastic origin: 
\begin{equation}
\label{equation "From the 4-vector acceleration"}
	\frac{1}{2}(DD^* + D^*D)\rotatedvarphi_t = (D_\textsc{s}D_\textsc{s} - D_\textsc{a}D_\textsc{a})\rotatedvarphi_t = D_\textsc{s}w^{\rotatedvarphi_t} - D_\textsc{a}w^{\rotatedvarphi_t},   
\end{equation}
constituting a Borel vector on a $\mathbb{R}$-field. From here, it is agile to have classical formulæ, 
\begin{align}
	D_\textsc{s}w^{\rotatedvarphi_t}_1 = \frac{\partial}{\partial_t}w^{\rotatedvarphi_t}_1 + \nabla_{w^{\rotatedvarphi_t}_1}w^{\rotatedvarphi_t}_1, \\ 
	D_\textsc{a}w^{\rotatedvarphi_t}_2 = \nabla_{w^{\rotatedvarphi_t}_2}w^{\rotatedvarphi_t}_2 + \frac{1}{2}\textcyrillic{\textit{в}}_\mathbbl{R}^2\nabla^2w^{\rotatedvarphi_t}_2. 
\end{align}
$\mathrm{N}\!\!\mathrm{B}$. Before arriving at an explicitly relativistic picture, the vector $w_1$ and $w_1$ can be written as 
\begin{align}
	w_1^{\rotatedvarphi_t} = w_1^{\rotatedvarphi_t}(m, t) = \frac{1}{2}\left[\rotatedPsi^0(m, t) + \rotatedPsi^0_*(m, t)\right], \\ 
	w_2^{\rotatedvarphi_t} = w_2^{\rotatedvarphi_t}(m, t) = \frac{1}{2}\left[\rotatedPsi^0(m, t) - \rotatedPsi^0_*(m, t)\right], 
\end{align}
respectively.

We can thereupon formalize the stochastic acceleration as follows:
\begin{equation}
\label{equation "Classical stochastic acceleration"}
	\frac{1}{2}(DD^* + D^*D)\rotatedvarphi_t = \left(\frac{\partial}{\partial_t}w^{\rotatedvarphi_t}_1 + \nabla_{w^{\rotatedvarphi_t}_1}w^{\rotatedvarphi_t}_1\right) - \left( \nabla_{w^{\rotatedvarphi_t}_2}w^{\rotatedvarphi_t}_2 + \frac{1}{2}\textcyrillic{\textit{в}}_\mathbbl{R}^2\nabla^2w^{\rotatedvarphi_t}_2\right).
\end{equation}

With the support of the latter equation, we are close to describing the accelerated motion $\ddot{x}(t)$ of a particle of mass $m$ and velocity $\dot{x}(t)$ along a curve $x(t)$, according to Newtonian mechanics: 
\begin{equation}
	\ddot{x} = 1/m\F\bigl(x(t), t, \dot{x}(t)\bigr),
\end{equation}
for a vector force field $\F$. If the system is conservative, Newton's second law, $\F = m\vec{a} = \frac{d\vec{v}}{dt} \equivalent \F\bigl(m(t), t, \dot{m}(t)\bigr) = \frac{D^{\mathring{\mathcal{T}}}}{dt}\dot{m}(t)$, becomes 
\begin{equation}
	\frac{D^{\mathring{\mathcal{T}}}}{dt}\dot{m}(t) = - \gradient{E_u}, 
\end{equation}	
where $D^{\mathring{\mathcal{T}}}/dt$ is the (covariant derivative of the) Levi-Civita connection, and $E_u$ is the potential energy.

From Eqq. \eqref{equation "From the 4-vector acceleration"} \eqref{equation "Classical stochastic acceleration"} we can extract this new equalization of the stochastic acceleration:\footnote{
	{} Compare with Eqq. \eqref{subequations "NMA"}. 
	}
\begin{equation}
	\frac{1}{2}(DD^* + D^*D)\rotatedvarphi_t = \left(\frac{\partial}{\partial_t}w_1^{\rotatedvarphi_t} + \nabla_{w_1^{\rotatedvarphi_t}}w_1^{\rotatedvarphi_t}\right) - \frac{1}{2}\left(2\nabla_{w_2^{\rotatedvarphi_t}}w_2^{\rotatedvarphi_t} - \Laplacian_{\textsc{rkh}}w_2\right),
\end{equation}
with the presence of the Laplace–de Rham operator, aka the Kodaira–Hodge Laplacian \cite[p. 196]{Dankel Jr. "Mechanics on Manifolds and the Incorporation of Spin into Nelson's Stochastic Mechanics"}, $\Laplacian_{\textsc{rkh}} = (d + \delta^2) = d\delta + \delta d$, including the Cartanian differential (namely, the exterior derivative) $d$, and the codifferential $\delta$.

\section{Relativistic Newton–Nelson Equations}
\label{section "Relativistic Newton–Nelson Equations"}

What we have seen in the previous Section leads to a double formula à la Nelson \cite{Nelson "Derivation of the Schrodinger Equation from Newtonian Mechanics"} \cite{Nelson "Construction of Quantum Fields from Markoff Fields"} \cite[chap. I, sec. 10. \textit{Stochastic Parallel Translation}]{Nelson "Quantum fluctuations"} \cite[chap. 12. \textit{Dynamics of stochastic motion}]{Nelson "Dynamical Theories of Brownian Motion"}, or à la Newton–Nelson,
\begin{align}
	\label{align "Newton–Nelson equations 1"}
	\frac{1}{2}\left(D^+D^- + D^- D^+\right)\rotatedvarphi_\tau & = \F\Bigl(\rotatedvarphi_\tau, \tilde{w}^\rotatedvarphi_1(\tau, \rotatedvarphi_\tau)\Bigr), \\
	\label{align "Newton–Nelson equations 2"}
	D^2\rotatedvarphi_\tau & = \frac{\hbar}{m}\idem_n, 
\end{align}
where

$\F$ is a linear operator $\F(x) \colon \mathcal{T}_x\mathcal{M}^4 \to \mathcal{T}_x\mathcal{M}^4$, physically interpretable as a 4-force (it is, again, a vector force field), imagining a stochastic mechanics in space-time of general relativity, with a relativistic particle endowed with rest mass $m$, under an Itô-like process $\rotatedvarphi_\tau$ in a 4-manifold,

$D^2\rotatedvarphi_\tau$ can be re-equalized (and clarified) like this,
\begin{equation}
	D^2\rotatedvarphi_\tau = \left\{\left[\lim_{\mathrm{\Delta}_\tau \downarrow 0}\textcyrillic{\textit{Э}}^\rotatedvarphi_\tau\left(\frac{\mathrm{\Delta}^+\rotatedvarphi_\tau \otimes \mathrm{\Delta}^+\rotatedvarphi_\tau}{\mathrm{\Delta}\tau}\right) = \lim_{\mathrm{\Delta}_\tau \downarrow 0}\textcyrillic{\textit{Э}}^\rotatedvarphi_\tau\left(\frac{\mathrm{\Delta}^-\rotatedvarphi_\tau \otimes \mathrm{\Delta}^-\rotatedvarphi_\tau}{\mathrm{\Delta}\tau}\right)\right] \mathrel{\bigg|}  \mathscr{U}^\rotatedvarphi_\tau\right\},
\end{equation}

$\idem_n$ is the identity matrix.

For the appearance of $\hbar$, see Eqq. \eqref{equation "Real constant (squared) and the h-bar"} \eqref{equation "Real constant and the h-bar"}. Consult T. Zastawniak \cite{Zastawniak "A Relativistic Version of Nelson's Stochastic Mechanics"}.

Two remarks.
\enumerationisinitium
\item In \eqref{align "Newton–Nelson equations 1"} \eqref{align "Newton–Nelson equations 2"} the \emph{Nelson's mean acceleration} (\textsc{nma}), which is a stochastic acceleration in the Nelsionan view, is shown; it can be established in various mathematical guises:
\begin{subequations}
\label{subequations "NMA"}
	\begin{empheq}[left = {\F^a_{\rotatedvarphi_\tau} = \empheqlbrace}]{align}
	& \tfrac{1}{2}(DD^* + D^*D){\rotatedvarphi_\tau}_{\left[\mathcal{M}^4\right]} = \tfrac{1}{2}\left(D^+D^- + D^- D^+\right){\rotatedvarphi_\tau}_{\left[\mathcal{M}^4\right]},\footnotemark \\
	& \tfrac{1}{2}\left[DD^*x(\tau) + D^*Dx(\tau)\right]\rotatedvarphi_\tau, \\
	& - \tfrac{1}{2}\left[\mathsf{B}_\tau\underdot{\mathsf{B}}_\tau(x) + \underdot{\mathsf{B}}_t\mathsf{B}_t(x)\right], \\
	& - \gradient{\OrnsteinUhlenbeck_\rotatedupsilon}/m,
    \end{empheq}
\end{subequations}
\footnotetext{
	{} Cf. \cite[p. 123]{Dohrn and Guerra "Nelson's stochastic mechanics on Riemannian manifolds"}.
		}
where 

$\mathsf{B}(t)$ and $\underdot{\mathsf{B}}(t)$ are differential operators inherent in Brownian motion, or Wiener process,

$\OrnsteinUhlenbeck_\rotatedupsilon$ is the potential under the Ornstein–Uhlenbeck process \cite{Uhlenbeck and Ornstein "On the Theory of the Brownian Motion"} addressed to the theory of Brownian motion.

It is very captivating to emphasize that Nelson \cite{Nelson "Derivation of the Schrodinger Equation from Newtonian Mechanics"} produces the Schrödinger equation 
\begin{subequations}
\label{subequations "Schrödinger equation"}	
\begin{align}
	i\hbar\frac{\partial\psi}{\partial t} & = i\frac{\textcyrillic{\textit{в}}^2_\mathbbl{R}}{2}\nabla^2\psi - i\frac{1}{\hbar}E_u\psi, \\
	& = -\frac{\hbar^2}{2}\nabla^2\psi + \underbrace{\gradient{E_u}}_{\gradient{\OrnsteinUhlenbeck_\rotatedupsilon}}\psi,\footnotemark \enspace \nabla^2 = \Laplacian \text{ (Laplace–Beltrami operator}),
\end{align}
\end{subequations}
\footnotetext{
	{} To be meticulous, $i\hbar\frac{\partial\psi}{\partial t}(x, t) = \text{ et cetera}$.
	}from diffusion theory, but he ends up developing a stochastic Newtonian equation. In other words, Eq. \eqref{subequations "Schrödinger equation"} is deductible from the stochastic Newton's equation. (And with that we reconnect to Section \ref{section "Some Tips. Let us Do a Refresh: Schrödinger Operators"}). This is verifiable because from a small modification of Eq. \eqref{equation "Classical stochastic acceleration"} and from 
\begin{align}
	\frac{1}{2}\left(DD^* + D^*D\right)\rotatedvarphi_t & = 1/m \cdot \F\Bigl(\rotatedvarphi_\tau, w^\rotatedvarphi_1(t, \rotatedvarphi_t)\Bigr), \\
	D^2\rotatedvarphi_t & = \frac{\hbar}{m}\idem_n,\footnotemark 
\end{align}
\footnotetext{
	{} Compare with Eqq. \eqref{align "Newton–Nelson equations 1"} \eqref{align "Newton–Nelson equations 2"}.
	}
revealing the trajectory of a particle under stochastic laws of motion, one sees that
\begin{equation}
	\frac{\partial w^\rotatedvarphi_1}{\partial t} = -\gradient{\OrnsteinUhlenbeck_\rotatedupsilon} -	\left(w^\rotatedvarphi_1\nabla\right)w^\rotatedvarphi_1 + \left(w^\rotatedvarphi_2\nabla\right)w^\rotatedvarphi_2 + \left(\frac{\textcyrillic{\textit{в}}^2_\mathbbl{R}}{2}\nabla^2\right)w^\rotatedvarphi_2,   
\end{equation}
for $\F = - \gradient{\OrnsteinUhlenbeck_\rotatedupsilon}$.
\item In \eqref{align "Newton–Nelson equations 1"} \eqref{align "Newton–Nelson equations 2"} we can also express the Ricci curvature, designated with $\Ric$, in the form of $\mathscr{C}^\infty$-smooth $\binom{1}{1}$-tensor,
\begin{align}
	\label{align "Newton–Nelson equations 3"}
	\frac{1}{2}(DD^* + D^*D)\rotatedvarphi_\tau & = 1/m \cdot \F\left(\rotatedvarphi_\tau, w^{\rotatedvarphi_\tau}_1, \tau\right) + \frac{\hbar}{2m}\Ric(\rotatedvarphi_\tau) \circ w_2^{\rotatedvarphi_\tau}, \\ 
	\label{align "Newton–Nelson equations 4"}
	D^2\rotatedvarphi_\tau & = \frac{\hbar}{m}g_{\{2, 0\}},
\end{align}
where $g_{\{2, 0\}}$ is an algebraic symmetric object in the form of $\binom{2}{0}$-tensor. Pay attention to the fact that, in the latter double formula, the mean derivatives and the Ricci curvature tensor have a definition according to the metric connection, viz. the Levi-Civita (Riemannian) connection 
\begin{equation}
	\omega_\mathfrak{h} \equivalent \nabla 
\end{equation}
of a $\binom{0}{2}$-tensor field, taking for granted that $\omega_\mathfrak{h} \in \Omega^1(\mathring{\mathcal{P}}, \mathfrak{h})$ is a $\mathfrak{h}$-valued 1-form on $\mathring{\mathcal{P}}$, or rather, a principal connection equivalent to the Levi-Civita connection (it is sometimes referred to as \emph{Ehresmann connection}).
\enumerationisfinis

\subsection{Parallel Translation of Random Vectors from One Fiber Bundle to Another}
\label{subsection "Parallel Translation of Random Vectors from One Fiber Bundle to Another"}

Let $(\mathring{\mathcal{E}}, \pi, \mathcal{M}, \mathring{\mathcal{F}})$ be a fiber bundle, where $\mathring{\mathcal{E}}$, $\mathcal{M}$ and $\mathring{\mathcal{F}}$ are the total space, the base space, and the fiber of the bundle, respectively. The vector bundle can be represented by the map $\pi \colon \mathring{\mathcal{E}} \to \mathcal{M}$. 

Presupposing a connection $\Gamma^{\mathring{\mathcal{E}}}_{\tau, s}$, the quartet of Eqq. \eqref{align "Newton–Nelson equations 1"} \eqref{align "Newton–Nelson equations 2"} \eqref{align "Newton–Nelson equations 3"} \eqref{align "Newton–Nelson equations 4"} can nimbly become a set of formulæ on fiber bundles. 

\vspace{2mm}

\begin{exemplum}
Let us confine ourselves to a paradigmatic equalities, with a 
\[
\text{$\left(\Gamma^{\mathring{\mathcal{E}}}_{\tau, s}\right)$-translation from } \mathring{\mathcal{E}}_{\rotatedvarphi_s} \text{ to } \mathring{\mathcal{E}}_{\rotatedvarphi_\tau}, 
\]
where the terms $\mathring{\mathcal{E}}_{\rotatedvarphi_s}$ and $\mathring{\mathcal{E}}_{\rotatedvarphi_\tau}$ are two fiber bundles related to stochastic processes with the occurrence of random variables. If we impose that 
\[
	\rotatedvarphi^{\mathring{\mathcal{E}}}_\tau \equivalent \rotatedvarphi_\tau 
\]
is our stochastic process in the bundle $\mathring{\mathcal{E}}$, then Eqq. \eqref{align "Forward mean derivative"} and \eqref{align "Backward mean derivative"}, articulating the covariant mean derivatives, take this look:
\begin{align}
	D^+\rotatedvarphi^{\mathring{\mathcal{E}}}_\tau = & \lim_{\mathrm{\Delta}_\tau \downarrow 0}\textcyrillic{\textit{Э}}\left\{\frac{\Gamma^{\mathring{\mathcal{E}}}_{\tau, \tau + \mathrm{\Delta}\tau}\rotatedvarphi^{\mathring{\mathcal{E}}}(\tau + \mathrm{\Delta}\tau) - \rotatedvarphi^{\mathring{\mathcal{E}}}_\tau}{\mathrm{\Delta}\tau} \mathrel{\bigg|} {\textit{\textgreek{Ν}}}^\rotatedvarphi_t\Bigl(\rotatedvarphi(\tau + \mathrm{\Delta}\tau) - \rotatedvarphi_\tau\Bigr)^2 \leqslant 0\right\} \notag \\
	& + \lim_{\mathrm{\Delta}_\tau\downarrow 0}\textcyrillic{\textit{Э}}\left\{\frac{\rotatedvarphi^{\mathring{\mathcal{E}}}_\tau - \Gamma^{\mathring{\mathcal{E}}}_{\tau, \tau - \mathrm{\Delta}\tau}\rotatedvarphi^{\mathring{\mathcal{E}}}(\tau - \mathrm{\Delta}\tau)}{\mathrm{\Delta}\tau} \mathrel{\bigg|}{\textit{\textgreek{Ν}}}^{\rotatedvarphi}_t\Bigl(\rotatedvarphi_\tau - \rotatedvarphi(\tau - \mathrm{\Delta}\tau)\Bigr)^2 \geqslant 0\right\}
\end{align}
and
\begin{align}
	D^-\rotatedvarphi^{\mathring{\mathcal{E}}}_\tau = & 
	\lim_{\mathrm{\Delta}_\tau \downarrow 0}\textcyrillic{\textit{Э}}\left\{\frac{\rotatedvarphi^{\mathring{\mathcal{E}}}_\tau - \Gamma^{\mathring{\mathcal{E}}}_{\tau, \tau - \mathrm{\Delta}\tau}\rotatedvarphi^{\mathring{\mathcal{E}}}(\tau - \mathrm{\Delta}\tau)}{\mathrm{\Delta}\tau} \mathrel{\bigg|} {\textit{\textgreek{Ν}}}^{\rotatedvarphi}_t\Bigl(\rotatedvarphi_\tau - \rotatedvarphi(\tau - \mathrm{\Delta}\tau)\Bigr)^2 \leqslant 0\right\} \notag \\
	& + \lim_{\mathrm{\Delta}_\tau \downarrow 0}\textcyrillic{\textit{Э}}\left\{\frac{\Gamma^{\mathring{\mathcal{E}}}_{\tau, \tau + \mathrm{\Delta}\tau}\rotatedvarphi^{\mathring{\mathcal{E}}}(\tau + \mathrm{\Delta}\tau) - \rotatedvarphi^{\mathring{\mathcal{E}}}_\tau}{\mathrm{\Delta}\tau} \mathrel{\bigg|} {\textit{\textgreek{Ν}}}^{\rotatedvarphi}_t\Bigl(\rotatedvarphi (\tau + \mathrm{\Delta}\tau) - \rotatedvarphi_\tau\Bigr)^2 \geqslant 0\right\}.
\end{align}
\end{exemplum}

\subsection{Stochastic Quantization of a Particle in a Non-Abelian Gauge Field}

In light of what has just been said in Section \ref{subsection "Parallel Translation of Random Vectors from One Fiber Bundle to Another"}, with an addition of the vector bundle $\mathring{\mathcal{E}}$, Eqq. \eqref{align "Newton–Nelson equations 1"} \eqref{align "Newton–Nelson equations 2"} are nothing more than a version of the \emph{quantum formulæ of the motion} of a certain particle \emph{in a gauge field}, under \emph{a stochastic process}. Here is how:
\begin{align}
	\label{align "Newton–Nelson equations 5"}
	\frac{1}{2}\left(D^{\mathring{\mathcal{E}}+}D^- + D^{\mathring{\mathcal{E}}-} D^+\right)\rotatedvarphi^{\mathring{\mathcal{E}}}_t & = \F_{t, \rotatedvarphi^{\mathring{\mathcal{E}}}_t}\left(w^{\Gamma(\vartheta, \Omega)}_{\rotatedvarphi^{\mathring{\mathcal{E}}}}\right), \\
	\label{align "Newton–Nelson equations 6"}
	D^2\rotatedvarphi_t & = \frac{\hbar}{m}\idem_n,
\end{align}
postulating a connection $\Gamma(\vartheta, \Omega)$ having 
 
a (sub)connection $\vartheta$-form, or a vector-valued 1-form, 
 
a curvature $\Omega$-form, or a 2-form of the Cartan connection, 

by virtue of which
 \begin{equation}
 	\Omega = D\vartheta, \enspace \Omega \equivalent \Omega^{{\mathring{\mathcal{E}}}}, \vartheta \equivalent \vartheta^{{\mathring{\mathcal{E}}}} 
 \end{equation}
 holds on $\mathring{\mathcal{E}}$. 
 
 $\mathrm{N}\!\!\mathrm{B}$. A kindred road-solution to this one is taken by Y.E. Gliklikh and N.V. Vinokurova \cite[p. 77]{Gliklikh and Vinokurova "The Newton-Nelson Equation on Fiber Bundles with Connections"}.

\section{Klein–Gordon \& Dirac Equations}

The quartet of Eqq. \eqref{align "Newton–Nelson equations 1"} \eqref{align "Newton–Nelson equations 2"} \eqref{align "Newton–Nelson equations 3"} \eqref{align "Newton–Nelson equations 4"} is plainly related to Klein–Gordon equation, which I will divide into three versions,
\begin{equation}
	\begin{rcases}
	\left(\frac{1}{c^2}\frac{\partial^2}{\partial{t^2}} - \nabla^2  + \frac{m^2c^2}{\hbar^2}\right)\psi \\
	\Bigl[\left(\omega_\mathfrak{h} - \frac{i}{\hbar}\tilde{w}, \omega_\mathfrak{h} - \frac{i}{\hbar}\tilde{w}\right) - \frac{1}{h^2}\Bigr]\psi \\
	\dAlembertian + m^2\psi
	\end{rcases}
 	= 0,
\end{equation}
where $\psi$ is a scalar/wave function, and $\dAlembertian = \frac{1}{c^2}\frac{\partial^2}{\partial{t^2}} - \nabla^2$ is the d'Alembertian, constituting the substratum for the Dirac equation, also divided into three versions,
\begin{equation}
\label{equation "Dirac equation"}
	\begin{rcases}
	\left(i\hbar\gamma^\mu\frac{\partial}{\partial{x}^\mu} - mc\right)\psi \\
	\bigl(i\gamma^\mu\partial_\mu - m\bigr)\psi \\
	\bigl(i\slashed{\partial} - m\bigr)\psi
	\end{rcases}
	= 0,\footnotemark
\end{equation}
where the Dirac 4-spinor, id est a 4-component wave function, $\psi = \left(\begin{smallmatrix}
	\zeta^\alpha \\ 
	\tilde{\chi}_{\dot{\alpha}}
	\end{smallmatrix}\right)$, appears, with a left-handed spinor $\zeta^\alpha$, and a right-handed spinor $\tilde{\chi}_{\dot{\alpha}}$.
\footnotetext{
	{} In Eq. \eqref{equation "Dirac equation"} $\gamma^\mu \text{ viz. } \gamma_\textsc{d}^\mu = \{\gamma^0, \gamma^1, \gamma^2, \gamma^3\}$ are the $[M]^{4 \times 4}$ Dirac gamma matrices, 
\begin{equation}
	\gamma^0 = 
	\begin{pmatrix*}[r]
	1 & 0 \\
	0 & -1
	\end{pmatrix*}, 
	\gamma^1 = 
	\begin{pmatrix}
	0 & \sigmaPauli_1 \\
	-\sigmaPauli_1 & 0
	\end{pmatrix},
	\gamma^2 = 
	\begin{pmatrix}
	0 & \sigmaPauli_2 \\
	-\sigmaPauli_2 & 0
	\end{pmatrix},
	\gamma^3 = 
	\begin{pmatrix}
	0 & \sigmaPauli_3 \\
	-\sigmaPauli_3 & 0
	\end{pmatrix},
	\end{equation}
plus 
$\gamma^5 = \bigl(\begin{smallmatrix}
	0 & 1 \\
	1 & 0
	\end{smallmatrix}\bigr)$; $\gamma^0$ is the time-like matrix, $1 = \idem_2$, and $\sigmaPauli_{1, 2, 3}$ the Pauli matrices; its extended writing is:
\begin{equation}
 	\gamma^0 = 
 	\Biggl\{\begin{smallmatrix} 
	1 & 0 & 0 & 0 \\
	0 & 1 & 0 & 0 \\ 
	0 & 0 & -1 & 0 \\
	0 & 0 & 0 & -1 
	\end{smallmatrix}\Biggr\},
	\gamma^1 =
	\Biggl\{\begin{smallmatrix}
	0 & 0 & 0 & 1 \\
	0 & 0 & 1 & 0 \\
	0 & -1 & 0 & 0 \\
	-1 & 0 & 0 & 0
	\end{smallmatrix}\Biggr\},
	\gamma^2 =
	\Biggl\{\begin{smallmatrix}
	0 & 0 & 0 & -i \\
	0 & 0 & i & 0 \\
	0 & i & 0 & 0 \\
	-i & 0 & 0 & 0
	\end{smallmatrix}\Biggr\}, 
	\gamma^3 =
	\Biggl\{\begin{smallmatrix}
	0 & 0 & 1 & 0 \\
	0 & 0 & 0 & -1 \\
	-1 & 0 & 0 & 0 \\
	0 & 1 & 0 & 0
	\end{smallmatrix}\Biggr\},
\end{equation} 
plus 
$\gamma^5 = i\gamma^0, \gamma^1, \gamma^2, \gamma^3 =
 	\Biggl\{\begin{smallmatrix} 
	0 & 0 & 1 & 0 \\
	0 & 0 & 0 & 1 \\ 
	1 & 0 & 0 & 0 \\
	0 & 1 & 0 & 0 
	\end{smallmatrix}\Biggr\}$, in the anti-commutation relation $\left\{\gamma^\mu, \gamma^\nu\right\} = \gamma^\mu\gamma^\nu + \gamma^\nu\gamma^\mu = 2\eta^{\mu\nu}\idem_{4 \times 4}$,

$\partial_\mu = \frac{\partial}{\partial{x}^\mu}$, and $\slashed{\partial} = \gamma^\mu\partial_\mu$ is the partial derivative under the Feynman slash notation, with which Feynman's Dirac operator, $D = \slashed{\partial} \equivalent \gamma^\mu\partial_\mu$, is accessed. 
	}

\section{Dirac Operator via Clifford Connection}

In these circumstances, we are able to find a link between the covariant derivatives of Sections \ref{subsection "Covariant Derivative of Stochastic Type (on a Flat Lorentz–Minkowski Space-Time) in Minkowski Space-Time via Itô–Wiener Processes"} plus \ref{subsection "Covariant Derivative of Stochastic Type in Curved Space-Time (on a Lorentz Hyperbolic Manifold) via Itô–Wiener Processes"} and the Dirac operator via Clifford algebra. Let

$\pi \colon \mathring{\mathcal{E}} \to \mathcal{M}$ be a vector bundle over a 4-manifold $\mathcal{M}$, 

$\Gamma_{\sezione}\bigl(\mathcal{M}, \mathring{\mathcal{E}}\bigr)$ be the space of smooth sections of $\mathring{\mathcal{E}}$,

$\Gamma_{\sezione}\bigl(\mathcal{M}, \bigwedge^k\mathring{\mathcal{T}}^*\mathcal{M} \otimes \mathring{\mathcal{E}}\bigr) \viz \textgreek{Ζ}_{\sezione}\bigl(\mathcal{M}, \mathring{\mathcal{E}}\bigr)$ be the space of sections of the bundle $\bigl(\mathcal{M}, \bigwedge^k\mathring{\mathcal{T}}^*\mathcal{M} \otimes \mathring{\mathcal{E}}\bigr)$, or the space of differential $k$-forms in $\mathring{\mathcal{E}}$, stated that $\mathring{\mathcal{T}}^*\mathcal{M}$ is a cotangent bundle (more particularly, a disjoint union of the cotangent spaces $\mathcal{T}^*_x\mathcal{M}$). 

So, assume that $\mathring{\mathcal{E}}$ is a complex $\mathbb{Z}_2$-bundle 
\begin{equation}
	\left\{\mathring{\mathcal{E}} \viz \bigl(\mathring{\mathcal{E}}\bigr)_{\mathbb{Z}_2}\right\} = \mathring{\mathcal{E}}^+ \oplus \mathring{\mathcal{E}}^-. 
\end{equation}
Which allows us to imagine that $\mathring{\mathcal{E}}$ is a bundle of Clifford modules, with a bundle map $\varphi_{\Cl} \colon \mathring{\mathcal{T}}^*\mathcal{M} \to \Endomorphism\bigl(\mathring{\mathcal{E}}\bigr)$ such that 
\begin{equation}
	\varphi_{\Cl}(\omega_1)\varphi_{\Cl}(\omega_2) + \varphi_{\Cl}(\omega_2)\varphi_{\Cl}(\omega_1) = -2(\omega_1, \omega_2), 
\end{equation}	
letting $\omega \in \mathcal{T}^*_x\mathcal{M}$ be a 1-form (acting as a cotangent vector) on $\mathcal{M}$.

Given a Cliffordian bundle of a Riemannian manifold $\Cl(\mathcal{M})$, and identified with $\Endomorphism_{\Cl(\mathcal{M})}$ an endomorphism of Clifford $\mathring{\mathcal{E}}$-valued module over $\Cl(\mathcal{M})$, the space of sections 
\begin{equation}
	\Gamma_{\sezione}\bigl\{\mathcal{M}, \Endomorphism\bigl(\mathring{\mathcal{E}}\bigr)\bigr\} \cong \Gamma_{\sezione}\bigl\{\mathcal{M}, \Cl(\mathcal{M}) \otimes \Endomorphism_{\Cl(\mathcal{M})}\bigl(\mathring{\mathcal{E}}\bigr)\bigr\} 
\end{equation}
is isomorphic to the space of $\mathring{\mathcal{E}}$-valued differential $k$-forms,
\begin{equation}
	\textgreek{Ζ}_{\sezione}\bigl\{\mathcal{M}, \Endomorphism_{\Cl(\mathcal{M})}\bigl(\mathring{\mathcal{E}}\bigr)\bigr\} \cong \Gamma_{\sezione}\left(\mathcal{M}, \bigwedge^k\mathring{\mathcal{T}}^*\mathcal{M} \otimes \Endomorphism_{\Cl(\mathcal{M})}\bigl(\mathring{\mathcal{E}}\bigr)\right).
\end{equation}

The endomorphic decomposition is 
\[
	\Endomorphism\bigl(\mathring{\mathcal{E}}\bigr) \cong \Cl(\mathcal{M}) \otimes \Endomorphism_{\Cl(\mathcal{M})}\bigl(\mathring{\mathcal{E}}\bigr).
\]

A connection on the vector bundle via Clifford $\mathring{\mathcal{E}}$-module, simply called \emph{Clifford connection}, is $\nabla^{\mathring{\mathcal{E}}}$, according to which 
\begin{equation}
	\left\{\nabla^{\mathring{\mathcal{E}}}_{\vec{X}}, \varphi_{\Cl}(\omega)\right\} = \varphi_{\Cl}\left(\nabla_{\vec{X}}\omega\right),
\end{equation}
where $\vec{X}$ is a vector field, and $\nabla_{\vec{X}}\omega$ is the Levi-Civita derivative of $\omega$. Which means that the Levi-Civita connection $\nabla^{\mathring{\mathcal{T}}}$ (to wit, the linear connection on the tangent bundle) on $\bigwedge^k\mathring{\mathcal{T}}^*\mathcal{M}$ is but a Clifford connection,
\begin{equation}
\label{equation "Levi-Civita and Clifford connection"}
	\nabla^{\mathring{\mathcal{T}}} \equivalent \nabla^{\mathring{\mathcal{E}}}.
\end{equation}

The \emph{Dirac operator} affiliated with the Clifford connection $\nabla^{\mathring{\mathcal{E}}}$ can take these two forms:
\enumerationisinitium
\item for a local frame field, namely an orthonormal basis, $e_j = \{e_1, \mathellipsis, e_n\}$, 
\begin{equation}
	D = \sum^n_{j = 1}\varphi_{\Cl}^j\nabla^{\mathring{\mathcal{E}}}_{e_j},
\end{equation}
or, alternatively, 
\begin{equation}
	D = \sum^n_{j = 1}\varphi_{\Cl}(dx^j)\nabla^{\mathring{\mathcal{E}}}_{\partial_j},
\end{equation}
\item for a sequential configuration,
\[
	\Gamma_{\sezione}\bigl(\mathcal{M}, \mathring{\mathcal{E}}\bigr) \xrightarrow{\nabla^{\mathring{\mathcal{E}}}} \Gamma_{\sezione}\bigl(\mathcal{M}, \mathring{\mathcal{T}}^*\mathcal{M} \otimes \mathring{\mathcal{E}}\bigr) \xrightarrow{\varphi_{\Cl}} \Gamma_{\sezione}\bigl(\mathcal{M}, \mathring{\mathcal{E}}\bigr).
\]
\enumerationisfinis

\vspace{2mm}

\begin{marginalia}[Smoothness, differentiability, and non-differentiability]
We were talking earlier on about \emph{space of sections}. We start from the surmise that our sections are \emph{smooth}. It should be remembered that by \emph{sections}—of vector bundles—we mean differentiable relations designed to colleague a certain vector in the corresponding vector space to each point of the manifold under consideration. This is in accordance with the scale relativity, which is a physico-geometric theory mixing differentiable and non-differentiable parts. The non-differentiable parts are certainly the \emph{fractal} ones.
\end{marginalia}

\vspace{2mm}

\begin{marginalia}[Spinor connection]
About Eq. \eqref{equation "Levi-Civita and Clifford connection"}, the same goes for a Levi-Civita connection $\nabla^{\mathring{\mathcal{P}}_\textit{ß}}$ on a spinor bundle $\mathring{\mathcal{P}}_\textit{ß}$ of a spin manifold $\mathcal{M}$. Here too we have a coincidence with the Clifford connection, 
\begin{equation}
\nabla^{\mathring{\mathcal{P}}_\textit{ß}} \equivalent \nabla^{\mathring{\mathcal{E}}}. 
\end{equation}
Recall that a spinor bundle $\mathring{\mathcal{P}}_\textit{ß}$ is 

a component of the principal $SL_2(\mathbb{C})$-bundle over $\mathbb{R}^4_{1, 3}$, or 

a complex vector bundle $\varsigma \colon \mathring{\mathcal{P}}_\textit{ß} \to \mathcal{M}$ (see Atiyah–Singer index theorem).
\end{marginalia}

\section{A Dive into the Fractoid Spaces}

In this Section, we will evaluate the notion of stochastic geodesic together with that of \emph{energy functional}, and of \emph{stochastic diffusions}, which is \emph{neither smooth nor time-differentiable}. (Since the geodesic paths, or the spatial trajectories of the manifold, are not differentiable anent the time, it is indispensable to substitute the \emph{derivative in time} with a \emph{mean derivative}). There will also be a nexus to Brownian motion. I call \emph{fractoid spaces} all non-time-differentiable fractal spaces that, wherever it happens, can also be devoid of smoothness.

\subsection{Derivatives Coupled with a I-Type Semi-martingale and the Stratonovich Integral}

Let $\textcyrillic{\textit{м}}_t$ be 

· a semi-martingale in $\mathcal{T}_{\textcyrillic{\textit{м}}_0}\mathcal{M}$, in accordance with a filtration $\filtration_t$, $t \in [0, T]$, and 

· a stochastic process with 
\begin{equation}
	\textcyrillic{\textit{м}}_t = \textcyrillic{\textit{м}}^\mathrm{loc}_t + \textgreek{\textit{χ}}_t,
\end{equation}	
where $\textcyrillic{\textit{м}}^\mathrm{loc}_t$ is a local martingale, and $\textgreek{\textit{χ}}_t$ is a stochastic, or random, process having paths almost surely bounded with variation in $t$. 

Let $\rotatedvarphi_t$ be a stochastic process outlined with a stochastic differential equation in Stratonovich form, and $\ell_s$ the arc length (the distance between two points along a section of a curve),
\begin{equation}
\label{equation "Stratonovich integral"}
	\rotatedvarphi_t = \int^t_0\Bigl(\tau^{\textcyrillic{\textit{м}}}_{0 \leftarrow \ell_s} \circ d\textcyrillic{\textit{м}}_{\ell_s}\Bigr),
\end{equation}
estimating that 
\[
	\tau^{\textcyrillic{\textit{м}}}_{t \leftarrow \ell_s} \colon \mathcal{T}_{\textcyrillic{\textit{м}}_{\ell_s}}\mathcal{M} \to \mathcal{T}_{\textcyrillic{\textit{м}}_t}\mathcal{M} 
\]
is the stochastic parallel transport along $\rotatedvarphi_t$ paired with the connection $\nabla^{\mathring{\mathcal{T}}}$. The \emph{time derivative} on the right-hand side of \eqref{equation "Stratonovich integral"} can be demarcated by \emph{conditional expectations under a change of probability measure} $\bbmu$ (on every Borel set, cf. Section \ref{subsection "Covariant Derivative of Stochastic Type (on a Flat Lorentz–Minkowski Space-Time) in Minkowski Space-Time via Itô–Wiener Processes"}),
\begin{equation}
\label{equation "Time derivative... conditional expectations"}
	\mathsf{D}_\bbmu\rotatedvarphi_t \equivalent D\rotatedvarphi_t = \lim_{\varepsilon \to 0}\frac{1}{\varepsilon}E^\rotatedvarphi_k\left(\rotatedvarphi(t + \varepsilon) - \rotatedvarphi_t \mathrel{\bigg|} \filtration_t\right), \enspace \rotatedvarphi_t = \rotatedvarphi(t + \varepsilon), 
\end{equation}
where $\mathsf{D}_\bbmu$ is the derivative with reference to probability $\bbmu$-measures, and $E^\rotatedvarphi_k$ is the stochastic kinetic energy of $\textcyrillic{\textit{м}}$. 

The derivative, in its \emph{forward version}, is, instead,
\begin{equation}
\label{equation "Derivative in its forward version"}
	\mathsf{D}^\nabla_\bbmu\textcyrillic{\textit{м}}_t \equivalent D^{\nabla^{\mathring{\mathcal{T}}}}\hspace{-4pt}\textcyrillic{\textit{м}}_t = \tau^{\textcyrillic{\textit{м}}}_{t \leftarrow 0}D\rotatedvarphi_t.
\end{equation}

To \emph{generalize} the derivative, merely subjoin a vector field, say, $\vec{X}$, 
\begin{equation}
\label{equation "generalized derivative"}
	\mathsf{D}^\nabla_\bbmu\vec{X}_t \equivalent D^{\nabla^{\mathring{\mathcal{T}}}}\hspace{-4pt}\vec{X}_t = \lim_{\varepsilon \to 0}{E^\rotatedvarphi_k}^{\filtration_t}\left\{\frac{\tau^{\textcyrillic{\textit{м}}}_{t \leftarrow t + \varepsilon}\vec{X}\bigl(t + \varepsilon, \textcyrillic{\textit{м}}(t + \varepsilon)\bigr) - \vec{X}(t, \textcyrillic{\textit{м}}_t)}{\varepsilon}\right\}.
\end{equation}

$\mathrm{N}\!\!\mathrm{B}$. Compare the above results with those in X. Chen and A.B. Cruzeiro \cite[sec. 2]{Chen Cruzeiro "Stochastic geodesics and forward-backward stochastic differential equations on Lie groups"}.

\subsection[The Concept of $\mathscr{E}$-Stochastic Geodesics and of $\upsilon^\rotatedvarphi$-Diffusion Process]{The Concept of $\protect\pseudobold{\mathscr{E}}$-Stochastic Geodesics and of $\protect\pseudobold{\upsilon^\rotatedvarphi}$-Diffusion Process}

\begin{figure}
\centering
	\includegraphics[width = 0.550\textwidth]{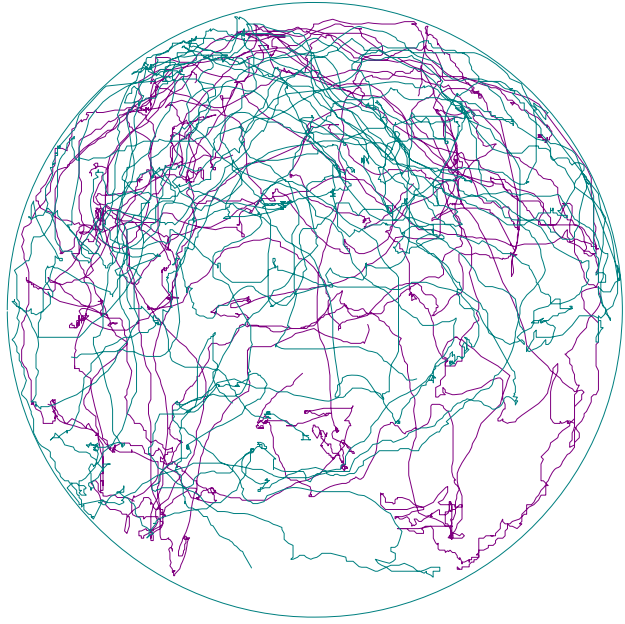}
	\caption{To give an idea with a simulation: Brownian-like \textcolor{eggplant}{\texttt{y-\#800080}} paths \& Brownian-like \textcolor{mallard}{\texttt{$\mathtt{x_t}$-\#008080}} paths in a circle, or a subset of the real 2-space, $\mathbb{S}^1 = \{x \in \mathbb{R}^2 \mid \|x\| = 1\}$}
	\label{figure "Brownian paths"} 
\end{figure}

The concept of geodesic can be characterized by the notion of set of critical points. The most direct manner is to use the \emph{random Brownian motion} (Fig. \ref{figure "Brownian paths"}), on the flat space, along certain directions, say, $\boundedvariation_\mathtt{s}$ ($\mathtt{s}$ is for $\mathtt{sign}$). Let $x$ reveal a path of a $\mathbb{R}$-valued Brownian motion, $x_t, t \in [0, T], x_0 = 0$. We can \emph{generalize} the derivative $\mathsf{D}_\bbmu \equivalent D$ on a flat space, which is that of random motion, to hold onto the stochastic schema.

\vspace{2mm}

\begin{propositio}[Geodesic as a set of critical points]
Let 
\begin{equation}
\label{equation "Energy functional with semi-martingale"}
\mathscr{E}^\rotatedvarphi_k(\textcyrillic{\textit{м}}) = E^\rotatedvarphi_k\int^T_0\left\|\mathsf{D}^\nabla_\bbmu\textcyrillic{\textit{м}}_t \equivalent D^{\nabla^{\mathring{\mathcal{T}}}}\hspace{-4pt}\textcyrillic{\textit{м}}_t\right\|^2dt 
\end{equation}
be the energy functional, $w$ a vector field, and $\upsilon^\rotatedvarphi$ the diffusion process (see Section \ref{subsubsection "Diffusion process and Brownian motion"}), produced by the operator $\textit{\L}_w$, $w \in \mathscr{C}^2$, see Eq. \eqref{equation "Lw operator plus smooth function"}. Iff 
\begin{equation}
\label{equations "equations = 0"}
	\begin{rcases}
	\nabla^{\mathring{\mathcal{T}}}w + \partial_tw + \frac{\Laplacian_w + \Ric_w}{2} \\
	D^{\nabla^{\mathring{\mathcal{T}}}}\hspace{-4pt}w(t, \upsilon^\rotatedvarphi_t)
	\end{rcases}
	= 0
\end{equation}
almost everywhere, said that $\Laplacian_w$ and $\Ric_w$\footnote{
	{} A definition of the Ricci curvature tensor for the stochastic geodesic curves inside the variational principle and the Euler–Lagrange dynamics, is in A.B. Cruzeiro \cite[p. 90]{Cruzeiro "Hydrodynamics Probability and the Geometry of the Diffeomorphisms Group"}.
	} 
are the Laplacian operator, see Eq. \eqref{equation "Laplacian on O(M)"}, and the Ricci curvature w.r.t. $w$, then $\upsilon^\rotatedvarphi$ is a critical path for $\mathscr{E}^\rotatedvarphi_k$, that is, a geodesic is nothing more than a set of critical points of 
\begin{equation}
\label{equation "Geodesic as a set of critical points"}
\mathscr{E}^\rotatedvarphi_k\left(\gamma_{\mathrm{c}(t)}\right) = \int^T_0\Bigl(g_{\mu\nu}\bigl(\dot{\gamma}_{\mathrm{c}(t)}\bigr)\dot{\gamma}_{\mathrm{c}(t)}^\mu\dot{\gamma}_{\mathrm{c}(t)}^\nu\Bigr)dt = \int^T_0\left\|\dot{\gamma}_{\mathrm{c}(t)}\right\|^2dt,
\end{equation}
$g_{\mu\nu}$ being the Riemannian metric tensor (field).
\end{propositio}
	
Eqq. \eqref{equation "Time derivative... conditional expectations"} \eqref{equation "Derivative in its forward version"} \eqref{equation "generalized derivative"} \eqref{equation "Energy functional with semi-martingale"} \eqref{equations "equations = 0"} \eqref{equation "Geodesic as a set of critical points"} can be placed into the \emph{fractoid spaces}.

\begin{proof}
Firstly, please note that the critical path satisfies the Euler–Lagrange equations:
	\begin{subequations}	
	\begin{align}
	& \frac{\partial\Lagrangian}{\partial x^\mu} = \frac{d}{dt}\left(\frac{\partial\Lagrangian}{\partial\dot{x}^\mu}\right) \text{ or } -\frac{\partial E_u}{\partial x^\mu} = \frac{d}{dt}\left(\frac{\partial E_k}{\partial\dot{x}^\mu}\right), \enspace 1 \leqslant \mu \leqslant n, \\
	& \frac{\partial\Lagrangian}{\partial x^\mu}\Bigl(\gamma_{\mathrm{c}(t)}, \dot{\gamma}_{\mathrm{c}(t)}, t\Bigr) = \frac{d}{dt}\left(\frac{\partial\Lagrangian}{\partial v^\mu}\right)\Bigl(\gamma_{\mathrm{c}(t)}, \dot{\gamma}_{\mathrm{c}(t)}, t\Bigr), \enspace \gamma_{\mathrm{c}(t)} = \left\{x^1_t, \mathellipsis, x^n_t\right\}, \dot{\gamma}_{\mathrm{c}(t)} = \left\{\dot{x}^1_t, \mathellipsis, \dot{x}^n_t\right\}.\footnotemark
	\end{align}
	\end{subequations}
\footnotetext{
	{} Consider a particle of mass $m$ moving along $\gamma_{\mathrm{c}(t)}$ with a velocity $v$. The function $\Lagrangian(x, \dot{x}, t) \text{ viz. } (x, v, t)$ is the Lagrangian on $\mathring{\mathcal{T}}\mathcal{M} \times [0, 1] \to \mathbb{R}$, and it matchs with the difference between the kinetic energy ($E_k$) and the potential energy ($E_u$).
	}
	
Let $\pi \colon \mathring{\mathcal{O}}(\mathcal{M}) \to \mathcal{M}$, $\pi(x, \textcyrillic{\textit{р}}) = x$, be a fiber map advantageous for a Euclidean isometry, and $\mathring{\mathcal{O}}(\mathcal{M}) = \left\{(x \in \mathcal{M}, \textcyrillic{\textit{р}})\right\}$ the orthonormal frame bundle on $\mathcal{M}$ (see Section \ref{subsubsection "Natural Isomorphism Allied with the Orthonormal Frame Bundle"}). If one handles the action functional of a Lagrangian system \cite[Eqq. (1.76) (1.77)]{Niccolai "Notes in Pure Mathematics and Mathematical Structures in Physics"}, and sets out that 
\begin{equation}
\label{equation "Energy functional"}
	\mathscr{E}^\rotatedvarphi_k = E^\rotatedvarphi_k\int^T_0\left\|D\pi \left(\textcyrillic{\textit{р}}_{x(t)}\right)\right\|^2dt,
\end{equation} 
one can obtain
\begin{align}
\label{align "Action functional..."}
	\frac{d}{d\varepsilon}\bigg|_{\varepsilon = 0}E^\rotatedvarphi_k\int^T_0\left\|D\pi \left(\textcyrillic{\textit{р}}_{x + \varepsilon\boundedvariation(t)}\right)\right\|^2dt & = 2E^\rotatedvarphi_k\int^T_0\left\{D\pi \left(\textcyrillic{\textit{р}}_{x(t)}\right), D\dot{\pi}\left(\frac{d}{d\varepsilon}\bigg|_{\varepsilon = 0}\left(\textcyrillic{\textit{р}}_{x + \varepsilon\boundedvariation(t)}\right)\right)\right\}dt \notag \\
	& = 2E^\rotatedvarphi_k\int^T_0\left\{D\pi \left(\textcyrillic{\textit{р}}_{x(t)}\right), \dot{\boundedvariation} - \frac{\Ric(\boundedvariation)}{2} - \Liederivative(\boundedvariation_t) \right\}dt \notag \\
	& = \Lbrack:2E^\rotatedvarphi_k\int^T_0\left\{D\pi \left(\textcyrillic{\textit{р}}_{x(t)}\right), D\left(\rotatedmcy_t\right)\right\}:\Rbrack dt,
\end{align}
where 

$\boundedvariation$ nominates processes of bounded variation, for which 
\begin{equation}
\label{equation "Processes of bounded variation"}
	\boundedvariation_0 = \boundedvariation_T = 0,
\end{equation}

$\Liederivative$ is the Lie derivative, cf. Eq. \eqref{equation "Laplacian on O(M)"}, 

the notations $\Lbrack:$ plus $:\Rbrack$ signify that the combination of symbols within them must be repeated, stimulated by the beginning and ending repeat signs in music, and

$\rotatedmcy$ is a semi-martingale, with initial conditions $\rotatedmcy_0 = 0$, see Eq. \eqref{equation "The semi-martingale rotatedmcy"}.

By virtue of the Eq. \eqref{equation "Processes of bounded variation"}, it follows that the second expression of \eqref{align "Action functional..."} $\Lbrack:\cdots:\Rbrack$ can be put in equivalence to $-2E^\rotatedvarphi_k\int^T_0\bigl\{(D)D\pi \left(\textcyrillic{\textit{р}}_x\right), \boundedvariation_t\bigr\} $, à savoir
\begin{equation}
	-2E^\rotatedvarphi_k\int^T_0\Bigl\{\left(D\right) D\pi \left(\textcyrillic{\textit{р}}_x\right), \boundedvariation_t\Bigr\} \equivalent 2E^\rotatedvarphi_k\int^T_0\left\{D\pi \left(\textcyrillic{\textit{р}}_{x(t)}\right), D\left(\rotatedmcy_t\right)\right\},
\end{equation}
and we are done. In Eqq. \eqref{equation "Energy functional"} \eqref{align "Action functional..."} \eqref{equation "Processes of bounded variation"} the equivalence $\mathsf{D}_\bbmu \equivalent D$ holds.
\end{proof}

\subsection{Marginalia of Clarification}

\subsubsection[The $\upsilon^\rotatedvarphi$-Diffusion Process and $\mathbb{R}$-valued Random Brownian Motion]{The $\protect\pseudobold{\upsilon^\rotatedvarphi}$-Diffusion Process and $\protect\pseudobold{\mathbb{R}}$-valued Random Brownian Motion}
\label{subsubsection "Diffusion process and Brownian motion"}

The diffusion process $\upsilon^\rotatedvarphi$ by $\textit{\L}_{w \in \mathscr{C}^2}$ is built on
\begin{equation}
	d{\upsilon^\rotatedvarphi_t}^{(\nu)} = y^\nu_\xi dx_t^\xi + dt\left\{w^\nu - \frac{g^{\varrho, \varsigma}\Gamma^\nu_{\varrho, \varsigma}}{2}\right\}; 
\end{equation} 
$y$ and $x_t$ are path $\mathbb{R}$-elements belonging to Brownian motion, accompanied by the Christoffel $\Gamma$-symbols.

Thinking back to the lesson of A.N. Kolmogorov \cite[p. 281]{"On the logarithmic normal distribution of particle sizes under grinding"}, a $\mathbb{R}$-valued random variable $y = \varphi(x)$ on a probability space $\mathscr{P}_\bbmu = (\invertedbreve{\Omega}, \mathscr{B}_\sigma, \bbmu)$ has a standard \emph{Gaussian distribution} with zero mean, where it takes the form
\begin{equation}
\label{equation "Gaussian distribution"}
	\bbmu\left(x \in \mathscr{P}_\bbmu \mid \varphi(x) < \rotatedeta\right) = \frac{1}{\sqrt{2\pi\left(\int_{\mathscr{P}_\bbmu}\varphi^2d\bbmu\right)}}\int^\rotatedeta_{-\infty}\exp\left\{-\frac{\textcyrillic{\textit{ш}}^2}{2\left(\int_{\mathscr{P}_\bbmu}\varphi^2d\bbmu\right)}\right\}d\textcyrillic{\textit{ш}}, \enspace -\infty < \rotatedeta < +\infty,
\end{equation}	
under a density function $\varphi(x)$.

The same formula applies to a space of (pseudo-)Riemannian type, replacing in \eqref{equation "Gaussian distribution"} $\mathscr{P}_\bbmu$ with the conventional \texttt{mathcal} letter $\mathcal{M}$, ou seja $x \in \mathcal{M}$, and $\int_\mathcal{M}$. This gives birth to what can be called \emph{$\mathbb{R}$-valued Brownian random measure}.
	
$\mathrm{N}\!\!\mathrm{B}$. From here it is possible to provide a subsequent definition of space decomposition of Hilbert space 
\begin{equation}
	\mathfrak{H} \equivalent \Lebesgue^2\left(\mathcal{M} \equivalent \mathbb{R}^n, \bbmu\right)
\end{equation}
into Hermite-Itô polynomial subspaces.

\subsubsection{Natural Isomorphism Allied with the Orthonormal Frame Bundle}
\label{subsubsection "Natural Isomorphism Allied with the Orthonormal Frame Bundle"}

Apropos of $\mathring{\mathcal{O}}(\mathcal{M})$, the orthonormal basis is $e(x)$ in $\mathcal{T}_x\mathcal{M}$, so that there is a natural isomorphism, in the Euclidean $n$-dimensional space, 
\begin{equation}
	\mathcal{T}_x\mathring{\mathcal{O}}(\mathcal{M}) \cong \mathbb{R}^n \oplus \mathfrak{o}_n(\mathbb{R}), 
\end{equation}
thru the medium of the Levi-Civita connection, and the orthogonal Lie algebra 
\begin{equation}
	\mathfrak{o}_n(\mathbb{R}) = \left\{X \in \mathfrak{gl}_n(\mathbb{R}) \mid X = -X^\textsc{t}\right\},\footnotemark
\end{equation}
for a skew symmetric matrix $X = -X^\textsc{t}$ in the Lie algebra $\mathfrak{gl}_n(\mathbb{R})$ of the general linear group $GL_n(\mathbb{R})$, set of all real $(n \times n)$-matrices.\footnotetext{
	{} $X^\textsc{t}$ is the transpose of $X$.
	}

\subsubsection{The II-Type Semi-martingale, and the Laplacian Operator}

The semi-martingale $\rotatedmcy$ is calculated by the method of Itô–Stratonovich, and it is definable as
\begin{equation}
\label{equation "The semi-martingale rotatedmcy"}
	\rotatedmcy = \frac{d}{d\varepsilon}\bigg|_{\varepsilon = 0}\left(\textcyrillic{\textit{р}}_{x + \varepsilon\boundedvariation}\right), \vartheta^\mathtt{l}_\mathtt{l},
\end{equation}	
once it is reported that $\vartheta^\mathtt{l}_\mathtt{l}$ ($\vartheta^\mathtt{l}$ or $\vartheta_\mathtt{l})$ is the habitual $\binom{1}{0}$-tensor valued 1-form ($\mathtt{l}$ is for $\mathtt{letter}$). Ergo (Section \ref{subsubsection "Natural Isomorphism Allied with the Orthonormal Frame Bundle"}) the Laplacian is so explicable, in relation to the Lie $\Liederivative$-derivative, 
\begin{equation}
\label{equation "Laplacian on O(M)"}
	\left(\Laplacian_w\right)_{\mathring{\mathcal{O}}(\mathcal{M})} = \sum^n_{\lambda = 1}\left(\Liederivative_{w_\vartheta}\right)^2.
\end{equation}
Implementing a $\mathscr{C}^\infty$ smooth function $\textcyrillic{\textit{г}}$, the Laplace–Beltrami operator on $\mathcal{M}$ is 
\begin{equation}
	\nabla^2 = \left(\Laplacian_w\right)_{\mathcal{M}}(\textcyrillic{\textit{г}})_{\mathscr{C}^\infty} = g^{\mu\nu}\left(\frac{\partial^2\textcyrillic{\textit{г}}}{\partial x^\mu\partial x^\nu} - \Gamma^\xi_{\mu, \nu}\frac{\partial\textcyrillic{\textit{г}}}{\partial x^\xi}\right), \enspace x \in \mathcal{M}.
\end{equation}

For a 4-dimensionality, in order to reconnect with Sections \ref{subsection "Covariant Derivative of Stochastic Type in Curved Space-Time (on a Lorentz Hyperbolic Manifold) via Itô–Wiener Processes"} and \ref{section "Relativistic Newton–Nelson Equations"}, we can straightforwardly write $\nabla^2 = \left(\Laplacian_w\right)_{\mathcal{M}^4}$.
	
And to complete the picture, the above-mentioned $\textit{\L}_w$ operator is so representable,
\begin{equation}
\label{equation "Lw operator plus smooth function"}
	\textit{\L}_w\textcyrillic{\textit{г}} = \frac{1}{2} \left(\Laplacian_w\right)_{\mathcal{M}} \textcyrillic{\textit{г}} + \partial_w\textcyrillic{\textit{г}}, \enspace w \in \mathscr{C}^2.
\end{equation}

\section[White Noise on a $(4+)\mathrm{D}$ Space-Time $\mathfrak{H}$-Geometry: the Paley–Wiener Integral]{White Noise on a $\protect\pseudobold{(4+)\mathrm{D}}$ Space-Time $\protect\pseudobold{\mathfrak{H}}$-Geometry: the Paley–Wiener Integral}

Section \ref{subsubsection "Diffusion process and Brownian motion"} sets the tone for bringing up Paley–Wiener's paraphernalia (stochastic calculus for Brownian motion) \cite{Paley Wiener and Zygmund "Notes on random functions"} \cite[chap. IX]{Paley and Wiener "Fourier Transforms in the Complex Domain"} \cite{Wiener "Nonlinear Problems in Random Theory"}. We may surmise that the Hilbert $\mathfrak{H}$-space has the following equivalence:
\begin{equation}
	\mathfrak{H} \equivalent \Lebesgue^2\left(\mathbb{R}^{4+}_{1, 3,}\right).\footnotemark
\end{equation}
\footnotetext{
	{} Pay attention to the comma after the number 3.
	}
So let 
\begin{equation}
	\mathbb{R}^{4+}_{1, 3,} = [0 , \infty) \times \mathbb{R}^{4+} 
\end{equation}
be a space-time, such as that in \cite{Niccolai "Spin and Torsion Tensors on Gauge Gravity: a Re-examination of the Einstein-Cartan Spatio-Temporal Theory"}, but still lacking the torsion, with a Minkowski–Lorentzian inner product $g(v, w) \text{ viz. } g\langle v, w\rangle$,\footnote{
	{} Isto é $g(v, w) = v^0w^0 + v^1w^1 + v^2w^2 + \cdots + v^{n - \mathbbl{z}_*}w^{n - \mathbbl{z}_*} - v^{n - \mathbbl{z}_* + 1}w^{n - \mathbbl{z}_* + 1} - \cdots - v^nw^n$, where $\mathbbl{z}_*$ is a non-negative integer.
	} 
\begin{equation}
	g(v, w)_{\langle v, w\rangle \in \mathfrak{H}} = \int^{+\infty}_0\int_{\mathbb{R}^{4+}_{1, 3,}}dtdx\left\{v_{t, x}w_{t, x}\right\}, \enspace t \geqslant 0, x \in \mathbb{R}^{4+}_{1, 3,}. 	
\end{equation}

Fix 
\enumerationisinitium
\item the process $\mathsf{W}_t = \{W_w\}_{w \in \mathfrak{H}}$\footnote{
	{} Be careful not to confuse this letter $W$ (via \texttt{\$W\$} command) with the letter $\Sobolev$ (via \texttt{\textbackslash{Sobolev}} command), for the Sobolev vector space, in Sec. \ref{subsection "The Head-scratcher of the Covariant Derivative"}: there is no correlation.
	} 
as a (continuous time) mean-zero gaussian distribution of the type of Eq. \eqref{equation "Gaussian distribution"}, having a \emph{covariance function} 
\begin{equation}
	\varphi_\textnormal{cov}(W_w, W_v) \equivalent g(w, v)_{\langle w, v\rangle \in \mathfrak{H}},
\end{equation}
\item and the white noise on $\mathbb{R}^{4+}_{1, 3,}$ as 
\begin{equation}
	\textgreek{\text{Θ}}^\text{w}_{t, x} = \frac{\partial^{4 + 1}W_{t, x}}{\partial_t\partial{x_1} \cdots \partial_{x_4}}. 
\end{equation}	
For a map $w \mapsto W_w$, the \emph{stochastic Paley–Wiener integration against} $\textgreek{\text{Θ}}^\text{w}$ is
\begin{equation}
	W_w = \int^{+\infty}_0\int_{\mathbb{R}^{4+}_{1, 3,}}dtdx\bigl\{w_{t, x}\textgreek{\text{Θ}}^\text{w}_{t, x}\bigr\}, \enspace w \in \mathfrak{H}.
\end{equation}
\enumerationisfinis

\section{Torsion of the Covariant Derivative}

To insert the spin-torsion tensor in \textsc{gr}, à savoir the Einstein–Cartan theory of gravity, see \cite{Niccolai "Spin and Torsion Tensors on Gauge Gravity: a Re-examination of the Einstein-Cartan Spatio-Temporal Theory"}. In detail, the Einstein–Cartan equations are those marked with the numbers (60) (61a) (61b) in \cite{Niccolai "Spin and Torsion Tensors on Gauge Gravity: a Re-examination of the Einstein-Cartan Spatio-Temporal Theory"}, whilst the Einstein–Cartan Lagrangian, for a $4\mathrm{D}$ space-time with spin-torsion (described by a Clifford $k$-form field) via Clifford bundles, is the one with the number (82).

Pertaining to the topic examined so far, we can write the \emph{torsion of the covariant derivative} in a variety of ways. I will highlight some of these.

\subsection[Torsion of $\nabla$-Connections thru Holonomy]{Torsion of $\protect\pseudobold{\nabla}$-Connections thru Holonomy}

We move on to the study of the holonomy of connections with skew symmetric torsion, as torsion is the skew symmetric contribution to the affine connection. Interestingly, the Einstein–Cartan theory can be deciphered as a theory of defects in a space-time with curvature and torsion, or in a 4-dimensional continuum containing defects. In fact, how much distance we keep, conceptually, between \emph{torsion} and \emph{deformation}? That is the question M.L. Ruggiero and A. Tartaglia \cite{Ruggiero and Tartaglia "Einstein-Cartan theory as a theory of defects in space-time"} are asking.

Let $\invertedbreve{\tau}$ be a torsion $\tbinom{1}{2}$-tensor of the $\nabla$-connection on a (pseudo-)Riemannian manifold, so we will take two $\mathscr{C}^\infty$ vector fields $\vec{X}$ and $\vec{Y}$, for which it holds that
\begin{equation}
	\invertedbreve{\tau}\left(\vec{X}, \vec{Y}\right) = \nabla_{\vec{X}}\vec{Y} - \nabla_{\vec{Y}}\vec{X} - \left[\vec{X}, \vec{Y}\right],
\end{equation} 
where

$\nabla_{\vec{X}}\vec{Y}$ is the covariant derivative of $\vec{Y}$ along $\vec{X}$, 

$\nabla_{\vec{Y}}\vec{X}$ is the covariant derivative of $\vec{X}$ along $\vec{Y}$, 

$\bigl[\vec{X}, \vec{Y}\bigr]$ is the commutator of $\vec{X}$ and $\vec{X}$.

Going to apply $\invertedbreve{\tau}$ to $\vec{X}$, one acquires a $\tbinom{1}{1}$-tensor $\invertedbreve{\tau}(\vec{X})$; if $\invertedbreve{\tau}(\vec{X})$ is applied to $\vec{Y}$, one acquires a $\tbinom{1}{0}$-tensor $\invertedbreve{\tau}(\vec{X}, \vec{Y})$.

In the presence of a scalar field a/o scalar function $\textgreek{\textit{\ddigamma}}^{\invertedbreve{\tau}}$, which is distinctly differentiable, the equivalence is
\begin{equation}
	\nabla_{\vec{X}}\left[\textgreek{\textit{\ddigamma}}^{\invertedbreve{\tau}}_{\vec{Y}}\right] = \vec{X}(\textgreek{\textit{\ddigamma}}^{\invertedbreve{\tau}})\vec{Y} + \textgreek{\textit{\ddigamma}}^{\invertedbreve{\tau}}\nabla_{\vec{X}}\vec{Y}. 
\end{equation}

\subsection[$D$-differentiation in Gravity Spin-Torsion Interaction thru Cartan–Einstein Model (in Riemann–Cartan Geometry)]{$\protect\pseudobold{D}$-differentiation in Gravity Spin-Torsion Interaction thru Cartan–Einstein Model (in Riemann–Cartan Geometry)}

Let 

$\textcyrillic{\textit{ч}}^4$ be a 4-form (which is the exterior derivative of a projective 3-form), 

$\Theta_{\invertedbreve{\tau}}$ be the torsion form, or the vector-valued 2-form, 

$\invertedbreve{\tau}$ be the torsion tensor (in this paper I use the notation \texttt{\textbackslash{tau}} plus an inverted breve, to avoid confusion with the notation \texttt{\textbackslash{tau}} connoting the proper time, see above).

Since 
\begin{equation}
	D\textcyrillic{\textit{ч}}^4_{\mu\nu\xi} = \textcyrillic{\textit{ч}}^4_{\mu\nu\xi\varrho}\Theta_{\invertedbreve{\tau}}^\varrho,
\end{equation}	 
adopting the algebraic Bianchi identities, one has initially these \emph{covariant exterior derivatives}
\begin{align}
	\label{align "torsion of D 1"}
	8\pi D_{\Tau_{\mu}} & = \frac{1}{2}\textcyrillic{\textit{ч}}^4_{\mu\nu\xi\varrho}\Theta_{\invertedbreve{\tau}}^\nu \wedge \Omega^{\xi\varrho}, \\
	\label{align "torsion of D 2"}
	8\pi D_{\hat{S}_{\mu\nu}} & = \textcyrillic{\textit{ч}}^4_{\nu\varrho} \wedge {\Omega^\varrho}_\mu - \textcyrillic{\textit{ч}}^4_{\mu\varrho} \wedge {\Omega^\varrho}_\nu,
\end{align}
where 

$\Tau^\mu = \Tau^{\mu\nu}\textcyrillic{\textit{ч}}^4_\nu$ is the tensor-valued 3-form, implying the energy-momentum tensor,

$\Omega^{\mathtt{l}\mathtt{l}}$, ${\Omega^{\mathtt{l}}}_{\mathtt{l}}$ is the curvature 2-form on space-time ($\mathtt{l}$ is for $\mathtt{letter}$, as already seen above), 

$\hat{S}_{\mu\nu}$ is the spin (density) tensor.

By selecting some suitable field equations,  such as
\begin{align}
	\frac{1}{2}g^{\xi\varrho}\textcyrillic{\textit{ч}}^4_{\mu\nu\xi} \wedge {\Omega^\nu}_\varrho & = -8\pi\Tau_{\mu}, \\
	\textcyrillic{\textit{ч}}^4_{\mu\nu\xi} \wedge \Omega^\xi & = 8\pi\hat{S}_{\mu\nu},
\end{align}
it is possible to draw the following rewriting of the previous Eqq. \eqref{align "torsion of D 1"} \eqref{align "torsion of D 2"},
\begin{align}
	D_{\Tau_\mu} & = {\invertedbreve{\tau}^\xi}_{\hspace{5pt}\mu\nu}\theta^\nu \wedge \Tau_\xi - \frac{1}{2}{\Riemann^\xi}_{\varrho\mu\nu}\vartheta^\nu \wedge {\hat{S}^\varrho}_\xi, \\
	D_{\hat{S}_{\mu\nu}} & = \vartheta_\nu \wedge \Tau_\mu - \vartheta_\mu \wedge \Tau_\nu, 
\end{align}
where $\theta^\mathtt{l}_\mathtt{l}$ is a 1-form, or a vector-valued 1-form.

\subsection{Quantum-like Fluctuations: a Stochastically Gravitational  Fabric of Space-Time} 

In this background it is subsequently permissible to include in our discussion
\enumerationisinitium
\item quantum-like fluctuations of the (pseudo-Euclidean) metric tensor $\eta_{\mu\nu}^{(1, 3)^+}$, $\eta_{\mu\nu}^{(1, 3)^-}$ of Minkowski space-time, and of the metric tensor $g_{\mu\nu}$ in Einstein's general relativity theory,
\item quantum-like fluctuations of the energy-momentum tensor $\Tau^{\mu\nu}$, $\Tau_{\mu\nu}$, oka stress-energy tensor, or stress-energy-momentum tensor,\footnote{
	{} $\Tau^{\mu\nu}$, or $\Tau_{\mu\nu}$, composes the kinetic energy of matter, since it acts as a matter-energy flow.
	} 
depicting the fluctuations of quantum matter fields in curved space-times,
\item quantum-like fluctuations of the whole gravitational field, taken in small (but how small?) pieces.
\enumerationisfinis

\vspace{2mm}

\addcontentsline{toc}{section}{Coda: Inspiring Snippet}
\begin{center}
\textbf{\textsc{coda: inspiring snippet}}
\end{center}

\vspace{2mm}

This paper has an essential bibliography as it is (was) a private communication, an \textit{échange de vues}, in order to be able to work with a certain lightness. Which results in a prehensility of the \textit{cerveau rêveur}. Rigor betwixt fantasy \& imagination. Is there a mathematical stream of consciousness? If so, this writing is a tiny epiphany of it.

There's a passage from A. Grothendieck \cite[6.2. (6). \textit{Le Rêveur}, p. 12 otm]{Grothendieck "Recoltes et Semailles. Reflexions et temoignage sur un passe de mathematicien"} that reads:

\vspace{2mm}

\begingroup
\footnotesize
Si nous pouvons communiquer avec nous-mêmes par le truchement du rêve, nous révélant à nous-mêmes, sûrement il doit être possible de façon toute aussi simple de communiquer à autrui le message nullement intime du rêve mathématique [\,\dots]. Et à vrai dire, qu'ai-je fait d'autre dans mon passé de mathématicien, si ce n'est suivre, “rêver” jusqu'au bout, jusqu'à leur manifestation la plus manifeste, la plus solide: irrécusable, des lambeaux de rêve se détachant un à un d'un lourd et dense tissu de brumes? Et combien de fois ai-je trépigné d'impatience devant ma propre obstination à polir jalousement jusqu'à sa dernière facette chaque pierre précieuse ou précieuse à demi en quoi se condensaient mes rêves — plutôt que de suivre une impulsion plus profonde: celle de suivre les arcanes multiformes du tissu-mère — aux confins indécis du rêve et de son incarnation patente, “publiable” en somme, suivant les canons en vigueur! J'étais d'ailleurs sur le point de suivre cette impulsion-là, de me lancer dans un travail de “science-fiction mathématique”, “une sorte de rêve éveillé”.

\endgroup

\setcounter{footnote}{0} 

\cleardoublepage
\phantomsection
\addcontentsline{toc}{section}{\refname}



\phantomsection
\addcontentsline{toc}{section}{Acknowledgements}

\thispagestyle{empty}

\begin{center}
\textsc{acknowledgements}
\end{center} 

\vspace{6mm}

\begingroup
\begin{center}
I owe a debt to my parents, Lucia and Bob. \\
I should like to thank my bruvver, Giulio, for his \textgreek{φιλάδελφος} impetus of generosity. \\
I am grateful also to my sis: Chiara had a hand in it, too. \\ 
There are, then, my consobrinus, Gabri, \& my senior confidantes, Filly, Laura, and Miranda, \\ 
who have shown themselves in full empathy with my woes. \\
The thoughtfulness of \ZhTraditional{埃琳娜} was revealed to be a sanative balsam. \\
All of them have sought to lend me their tootsies to carry on, \\ 
being my feet unusable. \\

\vspace{6mm}

\textcyrillic{Я благодарен \\ 
Эвелина Ю. Шамарова \\
за интерес, проявленный к этой математике}.	
\end{center}
\endgroup


\begin{thebibliography}{0}\small 
\thispagestyle{empty}
\markboth{Bibliography}{Bibliography}

\bibitem{Besov Il'In Nikol'skii "Integral Representations of Functions and Imbedding Theorems I"} O.V. Besov, V.P. Il'In, S.M. Nikol'skii, \textit{Integral Representations of Functions and Imbedding Theorems, Vol. I}, Ed. by M.H. Taibleson, J. Wiley \& Sons, Washington, 1978.
\bibitem{Besov Il'In Nikol'skii "Integral Representations of Functions and Imbedding Theorems II"} ------, \textit{Integral Representations of Functions and Imbedding Theorems, Vol. II}, Ed. by M.H. Taibleson, J. Wiley \& Sons, Washington, 1979.

\bibitem{Chen Cruzeiro "Stochastic geodesics and forward-backward stochastic differential equations on Lie groups"} X. Chen, A.B. Cruzeiro, \textit{Stochastic geodesics and forward-backward stochastic differential equations on Lie groups}, Discrete Contin. Dyn. Syst., Vol. 2013 (Supplement), pp. 115-121.

\bibitem{Cruzeiro "Hydrodynamics Probability and the Geometry of the Diffeomorphisms Group"} A.B. Cruzeiro, \textit{Hydrodynamics, Probability and the Geometry of the Diffeomorphisms Group}, in R.C. Dalang, M. Dozzi, F. Russo (Eds.), \textit{Seminar on Stochastic Analysis, Random Fields and Applications VI}, Centro Stefano Franscini, Ascona, May 2008, Birkhäuser · Springer Basel \textsc{ag}, Basel, 2011, pp. 83-93.

\bibitem{Cruzeiro and Zambrini "Feynman's Functional Calculus and Stochastic Calculus of Variations"} A.B. Cruzeiro and J.-C. Zambrini, \textit{Feynman's Functional Calculus and Stochastic Calculus of Variations}, in A.B. Cruzeiro, J.-C. Zambrini (Eds.), \textit{Stochastic Analysis and Applications}, Proceedings of the 1989 Lisbon Conference, Springer Science+Business Media, New York, 1991, pp. 82-95.

\bibitem{Dankel Jr. "Mechanics on Manifolds and the Incorporation of Spin into Nelson's Stochastic Mechanics"} T.G. Dankel Jr., \textit{Mechanics on Manifolds and the Incorporation of Spin into Nelson's Stochastic Mechanics}, Arch. Rational Mech. Anal., Vol. 37, № 3, 1970, pp. 192-221.

\bibitem{Dieudonne Schwartz "La dualite dans les espaces (F) et (LF)"} J. Dieudonné, L. Schwartz, \textit{La dualité dans les espaces $(\mathscr{F})$ et $(\mathscr{LF})$}, Ann. Inst. Fourier, Vol. 1, 1949, pp. 61-101.

\bibitem{Dohrn and Guerra "Nelson's stochastic mechanics on Riemannian manifolds"} D. Dohrn and F. Guerra, \textit{Nelson's stochastic mechanics on Riemannian manifolds}, Lett. Nuovo Cimento, Vol. 22, № 4, 1978, pp. 121-127.

\bibitem{Dohrn and Guerra "Geodesic correction to stochastic parallel displacement of tensors"} ------, \textit{Geodesic correction to stochastic parallel displacement of tensors}, in G. Casati, J. Ford (Eds.), \textit{Stochastic Behavior in Classical and Quantum Hamiltonian Systems}, Volta Memorial Conference, Como 1977, Springer-Verlag, Berlin, Heidelberg, 1979, pp. 241-249.

\bibitem{Dohrn Guerra Ruggiero "Spinning Particles and Relativistic Particles in the Framework of Nelson's Stochastic Mechanics"} D. Dohrn, F. Guerra, P. Ruggiero, \textit{Spinning Particles and Relativistic Particles in the Framework of Nelson's Stochastic Mechanics}, in S. Albeverio, Ph. Combe, R. Høegh-Krohn, G. Rideau, M. Sirugue-Collin, M. Sirugue and R. Stora (Eds.), \textit{Feynman Path Integrals}, Proceedings of the International Colloquium, Held in Marseille, May 1978, Springer-Verlag, Berlin, Heidelberg, 1979, pp. 165-181. 

\bibitem{Gliklikh and Vinokurova "The Newton-Nelson Equation on Fiber Bundles with Connections"} \textcyrillic{Ю.Е. Гликлих, Н.В. Винокурова}, \textcyrillic{\textit{Уравнение Ньютона–Нельсона на расслоениях со связностями}}, \textcyrillic{Фундаментальная и прикладная математика, том 20, № 3, 2015, pp. 61-81}.\footnote{ 
	{} An En. version, made by the authors themselves, Y.E. Gliklikh and N.V. Vinokurova, is in J. Math. Sci., Vol. 225, № 4, 2017, pp. 575-589, under the title \textit{The Newton–Nelson Equation on Fiber Bundles with Connections}.
	}

\bibitem{Grothendieck "Sur les espaces (F) et (DF)"} A. Grothendieck, \textit{Sur les espaces $(\mathscr{F})$ et $(\mathscr{DF})$}, Summa Bras. Math., Vol. 3, № 6, 1954, pp. 57-123. Fr. version is not found; version consulted: Ru. transl. by D.A. Raikov, Matematika, Vol. 2, № 3, 1958, pp. 81-127.
\bibitem{Grothendieck "Recoltes et Semailles. Reflexions et temoignage sur un passe de mathematicien"} ------, \textit{Récoltes et Semailles. Réflexions et témoignage sur un passé de mathématicien} [Juin 1983-Avril 1986], text in a free, online version: the number of pages refers to the page numbering of the manuscript (otm = of the manuscript).

\bibitem{Guerra and Ruggiero "A Note on Relativistic Markov Processes"} F. Guerra and P. Ruggiero, \textit{A Note on Relativistic Markov Processes}, Lett. Nuovo Cimento, Vol. 23, № 15, 9 Dic. 1978, pp. 529-534.

\bibitem{Ito "The Brownian motion and tensor fields on Riemannian manifold"} K. Itô, \textit{The Brownian motion and tensor fields on Riemannian manifold}, in  \textit{Proceedings of the International Congress of Mathematicians}, 15-22 August 1962, Stockholm, Institut Mittag-Leffler, Djursholm, Almqvist \& Wiksell, Uppsala, 1963, pp. 536-539.
\bibitem{Ito "Stochastic parallel displacement"} ------, \textit{Stochastic parallel displacement}, in M.A. Pinsky (Ed.), \textit{Probabilistic Methods in Differential Equations}, Proceedings of the Conference held at the University of Victoria, August 19-20, 1974, pp. 1-7.

\bibitem{"On the logarithmic normal distribution of particle sizes under grinding"} A.N. Kolmogorov, \textit{On the logarithmic normal distribution of particle sizes under grinding} (1941), in A.N. Shiryayev, \textit{Selected Works of A.N. Kolmogorov, Vol. II. Probability Theory and Mathematical Statistics}, Springer Science+Business Media, Dordrecht, 1992, pp. 281-284.

\bibitem{Littlewood Paley "Theorems on Fourier Series and Power Series"} J.E. Littlewood, R.E.A.C. Paley, \textit{Theorems on Fourier Series and Power Series}, J. Lond. Math. Soc. (1), Vol. 6, № 3, 1931, pp. 230-233.
\bibitem{Littlewood Paley "Theorems on Fourier Series and Power Series II"} ------, \textit{Theorems on Fourier Series and Power Series (II)}, J. Lond. Math. Soc. (2), Vol. 42, № 1, 1937, pp. 52-89.
\bibitem{Littlewood Paley "Theorems on Fourier Series and Power Series III"} ------, \textit{Theorems on Fourier Series and Power Series (III)}, J. Lond. Math. Soc. (2), Vol. 43, № 1, 1938, pp. 105-126.

\bibitem{Mandelbrot "Les objets fractals. Forme hasard et dimension"} B. Mandelbrot, \textit{Les objets fractals. Forme, hasard et dimension}, Flammarion, 1975, Paris (quatrième édition, 1995).
\bibitem{Mandelbrot "Measures of fractal lacunarity: Minkowski content and alternatives"} ------, \textit{Measures of fractal lacunarity: Minkowski content and alternatives}, in C. Bandt, S. Graf, M. Zähle (Eds.), \textit{Fractal Geometry and Stochastics}, Birkhäuser · Springer Basel \textsc{ag}, Basel, 1995, pp. 15-42.

\bibitem{Nelson "Derivation of the Schrodinger Equation from Newtonian Mechanics"} E. Nelson, \textit{Derivation of the Schrödinger Equation from Newtonian Mechanics}, Phys. Rev., Vol. 150, № 4, 1966, pp. 1079-1085. 
\bibitem{Nelson "Construction of Quantum Fields from Markoff Fields"} ------, \textit{Construction of Quantum Fields from Markoff Fields}, J. Funct. Anal., Vol. 12, № 1, 1973, pp. 97-112.
\bibitem{Nelson "Quantum fluctuations"} ------, \textit{Quantum fluctuations}, Princeton Univ. Press, Princeton (\textsc{nj}), 1985. 
\bibitem{Nelson "Dynamical Theories of Brownian Motion"} ------, \textit{Dynamical Theories of Brownian Motion}, re-typesetted second edition as \TeX{} file by J. Suzuki and revised by Nelson in 2001 (originally published by Princeton Univ. Press, 1967).

\bibitem{Niccolai "Notes in Pure Mathematics and Mathematical Structures in Physics"} E. Niccolai, \textit{Notes in Pure Mathematics \& Mathematical Structures in Physics}, \href{https://arxiv.org/abs/2105.14863}{arXiv:2105.14863} [math-ph], 2023 [v9]; the latest revision: \emph{Download} page at \href{https://www.edoardoniccolai.com}{https://www.edoardoniccolai.com}.
\bibitem{Niccolai "Spin and Torsion Tensors on Gauge Gravity: a Re-examination of the Einstein-Cartan Spatio-Temporal Theory"} ------, \textit{Spin \& Torsion Tensors on Gauge Gravity: a Re-examination of the Einstein–Cartan Spatio-Temporal Theory}, doi:\href{https://doi.org/10.5281/zenodo.7775360}{10.5281.7775360}, 2023 [v5], or \href{https://hal.science/hal-03948127}{hal-03948127} (HAL Id), 2023 [v3].

\bibitem{Nottale "Scale Relativity Fractal Space-Time and Quantum Mechanics"} L. Nottale, \textit{Scale Relativity, Fractal Space-Time and Quantum Mechanics}, Chaos Solitons Fractals, Vol. 4, № 3, 1994, pp. 361-388.
\bibitem{Nottale "Scale Relativity and Fractal Space-Time: Applications to Quantum Physics Cosmology and Chaotic Systems"} ------, \textit{Scale Relativity and Fractal Space-Time: Applications to Quantum Physics, Cosmology and Chaotic Systems}, Chaos Solitons Fractals, Vol. 7, № 6, 1996, pp. 877-938.
\bibitem{Nottale "Scale Relativity and Fractal Space-Time: Theory and Applications"} ------, \textit{Scale Relativity and Fractal Space-Time: Theory and Applications}, Found. Sci., Vol. 15, № 2, 2010, pp. 101-152; \href{https://arxiv.org/abs/0812.3857}{arXiv:0812.3857} [physics.gen-ph], 2008 [v1].
\bibitem{Nottale "Scale Relativity and Fractal Space-Time. A New Approach to Unifying Relativity and Quantum Mechanics"} ------, \textit{Scale Relativity and Fractal Space-Time. A New Approach to Unifying Relativity and Quantum Mechanics}, World Scientific, Singapore, 2011.

\bibitem{Nottale Celerier and T. Lehner Non-Abelian gauge field theory in scale relativity"} L. Nottale, M.-N. Célérier, and T. Lehner, \textit{Non-Abelian gauge field theory in scale relativity}, J. Math. Phys. 47, № 3, 2006, pp. 032303-1-19; \href{https://arxiv.org/abs/hep-th/0605280}{arXiv:hep-th/0605280}, 2006 [v1].

\bibitem{Nottale and Lehner "Turbulence and Scale Relativity"} L. Nottale and T. Lehner, \textit{Turbulence and Scale Relativity}, Phys. Fluids, Vol. 31, № 10, 2019, pp. 105109-1-22; \href{https://arxiv.org/abs/1807.11902}{arXiv:1807.11902} [physics.gen-ph], 2018 [v1].\bibitem{Ruggiero and Tartaglia "Einstein-Cartan theory as a theory of defects in space-time"} M.L. Ruggiero and A. Tartaglia, \textit{Einstein–Cartan theory as a theory of defects in space–time}, Amer. J. Phys., Vol. 71, № 12, 2003, pp. 1303-1313.

\bibitem{Paley and Wiener "Fourier Transforms in the Complex Domain"} R.E.A.C. Paley and N. Wiener, \textit{Fourier Transforms in the Complex Domain}, American Mathematical Society, New York, 1934.

\bibitem{Paley Wiener and Zygmund "Notes on random functions"} R.E.A.C. Paley, N. Wiener and A. Zygmund, \textit{Notes on random functions}, Math. Z., Vol. 37, № 1, 1933, pp. 647-668.

\bibitem{Sawano "Theory of Besov Spaces"} Y. Sawano, \textit{Theory of Besov Spaces}, Springer Nature, Singapore, 2018.

\bibitem{Schwartz "Produits tensoriels topologiques d'espaces vectoriels topologiques. Espaces vectoriels topologiques nucleaires"} L. Schwartz, \textit{Produits tensoriels topologiques d'espaces vectoriels topologiques. Espaces vectoriels topologiques nucléaires}, Séminaire Schwartz, Tome 1, 1953-1954, exp. № 1-24.

\bibitem{Sobolev "Sur un theoreme d'analyse fonctionnelle"} S.[L.] Sobolev (Soboleff), \textcyrillic{\textit{Об одной теореме функционального анализа}} (\textit{Sur un théorème d'analyse fonctionnelle}), Mat. Sb., Vol. 4(46), № 3, 1938, pp. 471-497.
\bibitem{Sobolev "Some Applications of Functional Analysis in Mathematical Physics"} ------, \textit{Some Applications of Functional Analysis in Mathematical Physics}, transl. from the Ru. by H.H. McFaden, American Mathematical Society (\textsc{ams}), Providence (\textsc{ri}), 2008\textsuperscript{re}.

\bibitem{Uhlenbeck and Ornstein "On the Theory of the Brownian Motion"} G.E. Uhlenbeck and L.S. Ornstein, \textit{On the Theory of the Brownian Motion}, Phys. Rev., Vol. 36, № 5, 1930, pp. 823-841.

\bibitem{Wiener "Nonlinear Problems in Random Theory"} N. Wiener, \textit{Nonlinear Problems in Random Theory}, The Technology Press of The Massachusetts Institute of Technology and J. Wiley \& Sons, New York, Chapman \& Hall, London, 1958.

\bibitem{Zastawniak "A Relativistic Version of Nelson's Stochastic Mechanics"} T. Zastawniak, \textit{A Relativistic Version of Nelson's Stochastic Mechanics}, Europhys. Lett., Vol. 13, № 1, 1990, pp. 13-17.
\end{thebibliography}
\end{document}